\documentclass{siamltex}
\usepackage{amsfonts}
\usepackage{amssymb,amsmath}
\usepackage{graphicx}
\usepackage{color}
\usepackage{hyperref}
\usepackage{subfigure}

\DeclareMathAlphabet{\itbf}{OML}{cmm}{b}{it}
\newcommand{\nc}{\newcommand}
\nc{\bx}{{\bf x}}
\nc{\by}{{\bf y}}
\nc{\bj}{{\bf j}}
\nc{\os}{\overline{\sigma}}
\nc{\de}{\delta}
\nc{\ep}{\varepsilon}

\nc{\cD}{\mathcal D}
\nc{\cU}{\mathcal U}
\nc{\bU}{\boldsymbol \cU}
\nc{\cE}{\mathcal E}
\nc{\cN}{\mathfrak{N}}
\nc{\cR}{\mathcal{R}}
\nc{\cX}{\mathcal{X}}
\nc{\cB}{\mathcal{B}}
\nc{\cT}{\mathcal{T}}
\nc{\bpsi}{\boldsymbol \Psi}

\begin{document}

\title{Asymptotic approximation of the Dirichlet to Neumann map of
  high contrast conductive media}

\author{ Liliana Borcea\footnotemark[1], Yuliya Gorb 
\footnotemark[2],
  and Yingpei Wang\footnotemark[3] }

\maketitle

\renewcommand{\thefootnote}{\fnsymbol{footnote}}

\footnotetext[1]{Computational and Applied Mathematics, Rice
  University, Houston, TX 77005. {\tt borcea@rice.edu}}
\footnotetext[2]{Department of Mathematics, University of Houston,
Houston, TX 77204. {\tt gorb@math.uh.edu}} 
\footnotetext[3]{Computational and Applied Mathematics, Rice
  University, Houston, TX 77005. {\tt Yingpei.Wang@rice.edu}}
\begin{abstract}
  We present an asymptotic study of the Dirichlet to Neumann map of
  high contrast composite media with perfectly conducting inclusions
  that are close to touching.  The result is an explicit
  characterization of the map in the asymptotic limit of the distance
  between the particles tending to zero.
\end{abstract}
\renewcommand{\thefootnote}{\arabic{footnote}}

\section{Introduction}
\label{sect:intro}
The Dirichlet to Neumann (DtN) map of an elliptic partial differential
equation maps the boundary trace of the solution to its normal
derivative at the boundary. It is used in inverse problems \cite{U-09}
for determining the coefficients of the elliptic equation, in non-overlapping domain decomposition methods \cite{TW-04} for solving
numerically the equations, and elsewhere. In this paper we study the
DtN map of equation
\begin{equation}
\nabla \cdot \left[ \sigma(\bx) \nabla u(\bx)\right] = 0, \quad \bx 
\in \cD\,, 
\label{eq:I1}
\end{equation}
with high contrast and rapidly varying nonnegative coefficient
$\sigma(\bx)$ in a bounded, simply connected domain $\cD \subset
\mathbb{R}^d$ with smooth boundary $\Gamma$. Rapidly varying means
that $\sigma$ fluctuates on a length scale that is much smaller than
the diameter of $\cD$. High contrast means that the ratio of the
largest and smallest value of $\sigma$ in $\cD$ is very large, even
infinite.  The coefficient $\sigma$ models the electrical conductivity
of a composite medium with highly conductive inclusions packed close
together in $\cD$, so that they are almost touching.  The solution $u$
of (\ref{eq:I1}) is the electric potential and $ -\sigma \nabla u$ is
the electric current, which we also call the flow.

The first mathematical studies of high contrast composites
\cite{Batchelor-77,Keller-63,Keller-87} are concerned with
homogenization of periodic media with perfectly conducting (or
insulating) inclusions.  Due to the periodicity, the problem reduces
to the local asymptotic analysis of the potential in the thin gap of
thickness $\delta$ between two neighboring inclusions. The asymptotics
is in the limit $\delta \to 0$.  The potential gradient in the gap
becomes singular in this limit, as described in
\cite{Ammari-05,Bao-09,Gorb-Novikov-12}, and the energy in the
composite is given to leading order by that in the gap, with effective
conductivity
\begin{equation}
\label{eq:I2}
\bar{\sigma} = \bar{\sigma}(\delta,g,d)\, .
\end{equation}
Here $d = 2$ or $3$ is the dimension of the space, and $g$ is a
geometrical factor depending on the local curvature of the boundaries
of the inclusions. The effective conductivity blows up in the limit
$\delta \to 0$ as $\delta^{-1/2}$ in two dimensions and
logarithmically in three dimensions.

Kozlov introduced in \cite{Kozlov-89} a continuum model of high
contrast conductivity in two dimensions
\begin{equation}
\label{eq:I3}
\sigma(\bx) = \sigma_o e^{S(\bx)/\epsilon}\, ,
\end{equation}
where $\sigma_o$ is a reference constant conductivity, $S(\bx)$ is a
smooth function with non-degenerate critical points, and $\epsilon \ll
1$ models the high contrast. An advantage of the model (\ref{eq:I3})
is that instead of specializing the analysis in the gaps to various
shapes of the inclusions, we can study a generic problem in the
vicinity of saddle points of the function $S(\bx)$.

In any case, independent of the model of high contrast, the problem
does not reduce to a local one if the medium does not have periodic
structure. The energy is still determined to the leading order by that in
the gaps, and each gap has an effective conductivity of the form
(\ref{eq:I2}), but the net flow in the gaps cannot be determined from
the local analysis.

The global problem is analyzed in \cite{BP-98}, for the high contrast
model (\ref{eq:I3}).  It uses two dual variational principles to
obtain sharp upper and lower bounds of the energy, which match to the
leading order.  The result can be interpreted as the energy of a
network with topology determined by the critical points of $S(\bx)$
i.e., $\sigma(\bx)$.  The nodes of the network are the maxima of
$S(\bx)$, and the edges connect the nodes through the saddle points of
$S(\bx)$.  Each saddle point $\bx_{_S}$ is associated with a resistor
with effective conductivity given by $\sigma(\bx_{_S})$ multiplied by
a geometrical factor depending on the curvatures of $S(\bx)$ at
$\bx_{_S}$.

The extension of the approach in \cite{BP-98} to homogenization of two
phase composites with infinite contrast is in \cite{Berlyand-01}. The
result is similar. The energy is given to leading order by a network
with nodes at the centers of the conductive inclusions. The edges
connect the nodes through the thin gaps separating the inclusions, and
have a net conductivity of the form (\ref{eq:I2}). An error analysis
of the approximation is in \cite{Novikov-02}.

The analysis of the DtN map is more involved than that of
homogenization, because of the arbitrary boundary conditions
\begin{equation}
u(\bx) = \psi(\bx), \quad \bx \in \Gamma.
\label{eq:I4}
\end{equation}
Still, the problems are related, because they both reduce to
approximating the energy in the composite, which can be bounded above
and below using dual variational principles. Indeed, the DtN map
$\Lambda:H^{1/2}(\Gamma) \to H^{-1/2}(\Gamma)$, defined by
\begin{equation}
\Lambda \psi (\bx) = \sigma(\bx) \nabla u(\bx) \cdot {\bf n}(\bx), 
\quad \bx \in \Gamma\, ,
\label{eq:I5}
\end{equation}
where ${\bf n}(\bx)$ is the outer normal at $\Gamma$, is self-adjoint.
Therefore, it is determined by its quadratic forms
\begin{equation}
  \left< \psi,\Lambda \psi \right> = \int_\Gamma ds(\bx) \, 
  \psi(\bx) \Lambda \psi(\bx)\, ,
\label{eq:I6p}
\end{equation}
for all $\psi \in H^{1/2}(\Gamma)$, and using integration by parts we
can relate it to the energy 
\begin{equation}
  E(\psi) = \frac{1}{2} \int_{\cD} d \bx \, \sigma(\bx) 
|\nabla u(\bx)|^2  \, ,
\label{eq:I7}
\end{equation}
by the equation
\begin{equation}
  \left< \psi,\Lambda \psi \right> = 2 E(\psi)\, .
\label{eq:I6}
\end{equation}

The DtN map of high contrast media with conductivity (\ref{eq:I3}) is
studied in \cite{BBP-99}. It is shown that $\Lambda$ can be
approximated by the matrix valued DtN map of the resistor network
described above, with topology determined by the critical points of
$S(\bx)$. However, the approximation in (\ref{eq:I3}) is on a subspace
of boundary potentials that vary slowly on $\Gamma$, on scales that
are larger or at most similar to the typical distance between the
critical points of $S(\bx)$.

In this paper we study the DtN map of two phase composites with
perfectly conducting inclusions in a medium of uniform conductivity
$\sigma_o = 1.$ For simplicity we work in two dimensions, in a disk
shaped domain $\cD$, with disk shaped inclusions.  The analysis
extends easily to any $\cD \subset \mathbb{R}^2$ with smooth boundary,
and to arbitrary inclusions, because only the curvature of their
boundary near the gaps plays a role in the approximation. The analysis
also extends to three dimensions, with some additional difficulties in
the construction of the test functions used in the variational
principles to bound the energy $E(\psi)$. High but finite contrast can
be handled as well, by writing the approximation as a perturbation
series in the contrast parameter, with terms calculated recursively,
as shown in \cite{calo-12,Gorb}.

As expected, we obtain that $\Lambda$ is determined by the DtN map of
the resistor network with nodes at the centers of the inclusions and
edges with effective conductivity of the form (\ref{eq:I2}). This is
the same network as in the homogenization studies
\cite{Berlyand-01,Novikov-02}. But the excitation of the network
depends on the boundary potential $\psi$. If $\psi$ varies slowly in
$\Gamma$, then the network plays the dominant role in the
approximation of $\Lambda$, and the result is similar to that in
\cite{BBP-99}. If $\psi$ varies rapidly in $\Gamma$, there is a
boundary layer of strong flow which must be coupled to the network.
The main result of the paper is the rigorous analysis of this
coupling. We show that the more oscillatory $\psi$ is, the less the
network gets excited, and the more dominant the boundary layer effect 
in the approximation of $\Lambda$. 

The paper is organized as follows: We begin in section \ref{sect:form}
with the formulation of the problem. The results are stated in section
\ref{sect:results} and the proofs are in section \ref{sect:proofs}. We
end with a summary in section \ref{sect:sum}.

\section{Formulation}
\label{sect:form}
\begin{figure}[t]
\begin{center}
\includegraphics[width=0.4\textwidth]{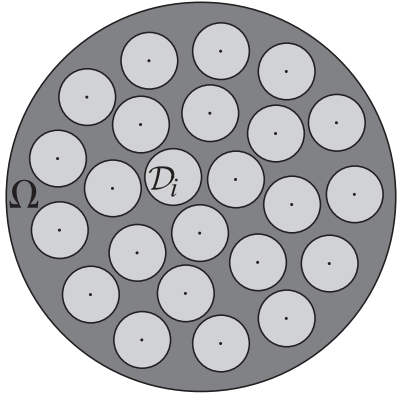}
\end{center}
\caption{Illustration of the setup. The domain $\cD$ contains $N$
  perfectly conducting inclusions denoted by $\cD_i$. The medium of
  conductivity $\sigma_o = 1$ lies in $\Omega$, the complement of the
  union of the inclusions in $\cD$.}
\label{fig:setup}
\end{figure}
We study the DtN map $\Lambda$ of an infinite contrast composite
medium in $\cD \subset \mathbb{R}^2$, consisting of $N \gg 1$
perfectly conducting inclusions $\cD_i$ centered at $\bx_i \in \cD$,
in a medium of uniform conductivity $\sigma_o$. See Figure
\ref{fig:setup} for an illustration. The domain $\cD$ is a disk of
radius $L$, centered at the origin of the system of coordinates. For
simplicity we let the inclusions be identical disks of radius $R \ll
L$. They are packed close together, but they are not touching. The
complement of the inclusions in $\cD$ is denoted by 
\[
\Omega = \cD 
  \setminus \bigcup_{i=1}^N \cD_i\, .
\]
\subsection{Variational principles}
\label{sect:VP}
The DtN map $\Lambda$ is determined by the quadratic forms
(\ref{eq:I6}), and therefore by the energy $E(\psi)$. We estimate it
using two dual variational principles. The first variational principle
\cite{BGN-05,Berlyand-01,Novikov-02}
\begin{equation}
  E(\psi) = \min_{v \in \mathbb{V}(\psi)} \frac{1}{2} \int_\Omega 
  d \bx 
  \, |\nabla v(\bx)|^2 \,, 
\label{eq:F1}
\end{equation}
is a minimization over potentials in the function space
\begin{equation}
  \mathbb{V}(\psi) = \left\{ v \in H^1(\Omega)\, , ~ v|_\Gamma = \psi\, 
    , ~     v|_{\partial \cD_i} = ~{\rm constant}, ~ ~i = 1, \ldots, N 
  \right\}\,.
\label{eq:F2}
\end{equation}
They have boundary trace $v|_\Gamma$ equal to the given $\psi \in
H^{1/2}(\Gamma)$, and are constant at the boundaries $\partial \cD_i$
of the inclusions. There is a unique minimizer of (\ref{eq:F1}), the
solution of the Euler-Lagrange equations \cite{Berlyand-01}
\begin{align}
  \Delta u(\bx) &=0\, , \qquad ~\bx \in \Omega\, , \label{eq:F3} \\
  u(\bx) &= \cU_i\,,\quad ~~ \bx \in \partial \cD_i\, , \label{eq:F4} \\
  \int_{\partial \cD_i} ds(\bx) \, {\bf n}(\bx) \cdot \nabla u(\bx)
  &=0\,, \qquad ~i = 1, \ldots, N\,, \label{eq:F5}\\
  u(\bx) &= \psi(\bx)\, , ~ ~\bx \in \Gamma \,.
\label{eq:F6}
\end{align}
The unknowns in these equations are the potential function $u(\bx)$
and the vector $\bU = (\cU_1, \ldots, \cU_N)$ of constant potentials
on the inclusions. These are the Lagrange multipliers associated with
the conservation of current conditions (\ref{eq:F5}).

The second variational principle
\begin{equation}
  E(\psi) = \max_{\bj \in \mathbb{J}} \left[\int_\Gamma ds(\bx) \, 
    \psi(\bx) {\bf n}(\bx) \cdot \bj(\bx) - \frac{1}{2}\int_\Omega 
    d\bx \, \left|\bj(\bx)\right|^2 \right]\, ,
\label{eq:F7}
\end{equation}
is a maximization over fluxes ${\bf j}$ in the function space
\begin{equation}
  \mathbb{J} = \left\{ \bj \in L_2(\Omega)\,, ~ \nabla \cdot \bj =0 
    ~ ~ {\rm in} ~ \Omega\, , 
    ~ \int_{\partial \cD_i} ds \, {\bf n} \cdot \bj = 0\, , ~ ~ i = 1, 
    \ldots, N \right\}\, .
\label{eq:F8}
\end{equation}
It is obtained from (\ref{eq:F1}) using Legendre (duality)
transformations \cite{Ekeland-76}, as explained for example in
\cite{Berlyand-01}.  The divergence free condition on $\bj$,
interpreted in the weak sense, gives the conservation of current in
$\Omega$, and the constraints at $\partial \cD_i$ are the analogues of
(\ref{eq:F5}). There is a unique maximizer of (\ref{eq:F7}), given by
\begin{equation}
\bj(\bx) = \nabla u(\bx),
\label{eq:F9}
\end{equation}
in terms of the solution of (\ref{eq:F3})-(\ref{eq:F6}).  It is the
negative of the electric current in $\Omega$. 

If we could solve equations (\ref{eq:F3})-(\ref{eq:F6}), we would have
the exact energy. This is impossible analytically.  Moreover,
numerical approximations of $(u,\bU)$ are computationally intensive
due to fine meshes needed to resolve the flow between the inclusions,
and the poor condition numbers of the resulting linear systems.  We
use instead the variational principles (\ref{eq:F1}) and (\ref{eq:F7})
with carefully constructed test functions $v \in \mathbb{V}(\psi)$ and
${\bf j} \in \mathbb{J}$ to obtain tight upper and lower bounds on
$E(\psi)$, which match to leading order.  The test potentials $v$ are
pieced together from local approximations of the solution of
(\ref{eq:F3})-(\ref{eq:F6}) in the gaps between the inclusions and in
a boundary layer at $\Gamma$. The construction of the test fluxes is
based on the relation (\ref{eq:F9}) between the optimal potential and
flux. Once we have a good test potential $v \in \mathbb{V}$, we can
construct ${\bf j} \in \mathbb{J}$ so that ${\bf j} \approx \nabla v$.

\subsection{Asymptotic scaling regime}
\label{sect:form.Sc}
There are three important length scales in the problem: The radius $L$
of the domain $\cD$, the radii $R$ of the inclusions and the typical
distance $\delta$ between the inclusions. To define $\delta$, we
specify first what it means for two inclusions to be neighbors.

Let $\mathfrak{X}_i$ be the Voronoi cell associated to the $i-$th
inclusion
\[
\mathfrak{X}_i = \left\{ \bx \in \cD ~\mbox{such that} ~ |\bx-\bx_i| \le
  |\bx-\bx_j|\, ,~~ \forall j = 1, \ldots, N, ~ j \ne i \right\}\,.
\]
It is a convex polygon bounded by straight line segments called edges.
The inclusion $\cD_i$ neighbors $\cD_j$ if the cells ${\mathfrak{X}_i}$
and ${\mathfrak{X}_j}$ share an edge.  We denote the set of indices of
the neighbors of $\cD_i$ by $\mathfrak{N}_i$,
\begin{equation}
\label{eq:F10}
\mathfrak{N}_i = \left\{ j\in \{1, \ldots, N\}\, , \quad
 j ~{\rm neighbors} ~ i \right\} \, ,
\end{equation}
and let 
\begin{equation}
 \label{eq:F11}
\delta_{ij} = \mbox{dist}\{\cD_i, \cD_j\}\, ,
\end{equation}
for all $i = 1, \ldots, N$ and $j \in \mathfrak{N}_i$. These are the
thicknesses of the gaps between the inclusions.

Similarly, we say that inclusion $\cD_i$ neighbors the boundary 
if ${\mathfrak{X}_i} \cap \Gamma \ne \emptyset$. Let us say that 
there are $N^\Gamma$ such inclusions and let $\delta_i$
be their distance from the boundary
\begin{equation}
\delta_i = \mbox{dist}\{\cD_i,\Gamma\}\,.
\label{eq:F12}
\end{equation}
We number henceforth the inclusions starting with those neighboring
$\Gamma$, counterclockwise. Thus, $\cD_i$ neighbors $\Gamma$ if $i =
1, \ldots, N^\Gamma$, and it is an interior inclusion if $i =
N^\Gamma+1, \ldots, N$.

We assume that both $\delta_{ij}$ and $\delta_i$ are of the same order
$\delta$, and seek an approximation of the DtN map $\Lambda$ in the
asymptotic regime of separation of scales
\begin{equation}
\delta \ll R \ll L.
\label{eq:F13}
\end{equation}
The reference order one scale is $L$. 

There is one more parameter in the asymptotic analysis, denoted by
$k$, which defines the Fourier frequency of oscillation of $\psi$ at
$\Gamma$.  It is independent of all the other scales in the problem
and it can vary between $0$ and $K$, with $K$ arbitrarily large. For
example, in domain decomposition, $K$ would be determined by the mesh
used to discretize the domain.  Because $\Gamma$ is a circle of radius
$L$ in our setup, we parametrize it by the angle $\theta \in [0,
2\pi]$, and suppose that $\psi$ is a superposition of Fourier modes
\begin{equation}
\psi(\theta) = \sum_{k=0}^K \left[ a_k^c \cos(k \theta) + 
a_k^s \sin(k \theta) \right]\, .
\label{eq:FDec}
\end{equation}
We seek approximations of $\left< \psi,\Lambda \psi\right>$ that are
valid for any $K$.

\section{Results}
\label{sect:results}
We state in Theorem \ref{thm.1} the approximation of
$\left<\psi,\Lambda \psi \right>$ for boundary potentials $\psi$ given
by a single Fourier mode.  The generalization to potentials
(\ref{eq:FDec}) is in Corollary \ref{cor.1}

The approximation involves the discrete energy and therefore DtN map
of a resistor network that is uniquely determined by the medium.  It
has the graph $\left(\cX,\mathfrak{E}\right)$ and edge conductivity
function $\os: \mathfrak{E} \to \mathbb{R}^+$. Each edge is associated
to a gap between adjacent inclusions or between an inclusion and the
boundary, and models the net singular flow there. The set of nodes of
the network is given by
\begin{equation}
\cX = \left\{ \bx_i, \quad i = 1, \ldots, N, ~ ~ 
\bx_i^\Gamma, \quad i = 1, \ldots, N^\Gamma\right\}\, .
\label{eq:R5}
\end{equation}
The interior nodes $\bx_i$ are at the centers of the inclusions, for
$i = 1, \ldots, N$. The boundary nodes
\begin{equation}
\bx_i^\Gamma = L(\cos \theta_i,\sin \theta_i) \,
\label{eq:R6}
\end{equation}
are the closest points on $\Gamma$ to the inclusions $\cD_i$ in its
vicinity, for $i = 1, \ldots, N^\Gamma$. The edges of the network
connect the adjacent nodes
\begin{equation}
  \mathfrak{E} = \left\{ e_{ij} = (\bx_i, \bx_j), ~ ~ i = 1, \ldots, 
    N, ~ ~ j \in \cN_i, \quad e_i^\Gamma = (\bx_i^\Gamma,\bx_i), ~ ~ 
i = 1, \ldots N^\Gamma \right \} \, ,
\label{eq:R7}
\end{equation}
and the network conductivity function is defined by
\begin{align}
   \os ( e_{ij} ) &= \pi \sqrt{ \frac{R}{\de_{ij}}} =: \os_{ij}, \quad i
  = 1,
  \ldots, N, ~ ~ j \in \cN_i\, , \label{eq:R8}\\
   \os(e_{i}^\Gamma) &= \pi \sqrt{\frac{2 R}{\de_{i}}} =: \os_i, \quad
  i = 1, \ldots, N^\Gamma \, .\label{eq:R9}
\end{align}

The DtN map $\Lambda^{\rm net}$ of the network is a symmetric
$N^\Gamma \times N^\Gamma$ matrix. Its quadratic forms are related to
the discrete energy $E^{\rm net}(\bpsi)$ of the network by
\begin{equation}
  \bpsi \cdot \Lambda^{\rm net} \bpsi = 2 E^{\rm net}(\bpsi)\, ,
\label{eq:RN1}
\end{equation}
where we let $\bpsi = (\Psi_1, \ldots,
\Psi_{N^\Gamma})^T$ be the vector of boundary potentials. The
energy has the variational formulation
\begin{equation}
  E^{\rm net}(\bpsi) = 
  \min_{\bU \in \mathbb{R}^{N}} \left\{
    \sum_{i=1}^{N^\Gamma} \frac{\os_i}{2} \left[ \cU_i -
      \Psi_i \right]^2 + \frac{1}{2}\sum_{i=1}^N
    \sum_{j\in \cN_{i}} \frac{\os_{ij}}{2} 
    (\cU_i-\cU_j)^2 \right\}\, ,
\label{eq:R4}
\end{equation}
where the $1/2$ factor in front of the second sum is because we sum
twice over the edges $e_{ij}$.  There is a unique minimizer $\bU \in
\mathbb{R}^{N}$ of (\ref{eq:R4}). It is the vector of node potentials
that satisfy Kirchhoff's equations, a linear system which states that
the sum of currents in each interior node equals zero.

\subsection{Boundary potential given by a single Fourier mode}
\label{sect:SFM}
Let the boundary potential $\psi$ be given by 
\begin{equation}
\psi(\theta) = \cos(k \theta) \, ,
\label{eq:R1}
\end{equation}
with $k > 0$. The case $k = 0$ is trivial, because constant
potentials are in the null space of the DtN map. 
\begin{theorem}
\label{thm.1}
We have that
\begin{equation}
  \left<\psi,\Lambda \psi \right> = 2 E(\psi) = 2 \cE(\psi)
  \left[ 1 +    o(1)\right]\, ,
\label{eq:R2}
\end{equation}
with the leading order of the energy given by the sum of three terms
\begin{equation}
\label{eq:R3}
\cE(\psi) = E^{\rm net}\left(\bpsi(\psi)\right) + 
\frac{k \pi}{2} + \cR_k \, . 
\end{equation}
The first term is the discrete energy $E^{\rm net}(\bpsi(\psi))$ of
the resistor network described above in \eqref{eq:R4}, with vector $\bpsi =
\left(\Psi_1, \ldots, \Psi_{N^\Gamma}\right)^T$ of boundary potentials
defined by 
\begin{equation}
  \Psi_i(\psi) = \psi(\theta_i) e^{-\frac{k \sqrt{2 R \de_i}}{L}},
  \quad i = 1, \ldots, N^\Gamma.
\label{eq:RN3}
\end{equation}
The second term in (\ref{eq:R3}) is the energy in the reference medium
with constant conductivity $\sigma_o = 1$. It is related to the
reference DtN map $\Lambda_o$ by
\begin{equation}
E_o(\psi) = \frac{1}{2} \int_\cD d \bx \, |\nabla u_o(\bx)|^2 = 
\frac{k \pi}{2} = \frac{1}{2} \left< \psi,\Lambda_o \psi\right>\, . 
\label{eq:R10}
\end{equation}
The last term in (\ref{eq:R11}) is given by
\begin{equation}
  \cR_k = \sum_{i=1}^{N^\Gamma}  \frac{\os_i}{4} 
  \left[
    \sqrt{\frac{2 k \de_i}{\pi L}} {\rm Li}_{1/2} 
    \left(e^{-\frac{2 k \de_i}{L}}\right) - e^{-\frac{2k 
        \sqrt{2 R
          \de_i}}{L}}\right]\, ,
\label{eq:R11}
\end{equation}
in terms of the Polylogarithm function ${\rm Li}_{1/2}$.
\end{theorem}

The proof of the theorem is in section \ref{sect:proofs}, and the
meaning of the result is as follows. The resistor network plays a role
in the approximation if it gets excited. This happens when the
boundary potential $\psi$ is not too oscillatory. As shown in equation
(\ref{eq:RN3}), the potential $\Psi_i$ at the $i-$th boundary node of
the network is not simply $\psi(\theta_i)$. We have an exponential
damping factor, which is due to the fact that only part of the flow
reaches the inclusion $\cD_i$. As $k$ increases, the flow near the
boundary becomes oscillatory, and has a strong tangential component.
Less and less current flows into $\cD_i$ and in the end, the network
may not even get excited. 

The term $\cR_k$ in (\ref{eq:R3}), which we rewrite as 
\begin{equation}
\cR_k = \sum_{i=1}^{N^\Gamma} \cR_{i,k}\, ,
\end{equation}
with 
\begin{equation}
\label{eq:defRik}
\cR_{i,k} = \frac{\os_i}{4} 
  \left[
    \sqrt{\frac{2 k \de_i}{\pi L}} {\rm Li}_{1/2} 
    \left(e^{-\frac{2 k \de_i}{L}}\right) - e^{-\frac{2k 
        \sqrt{2 R
          \de_i}}{L}}\right]\, , 
\end{equation}
describes the anomalous energy due to the oscillations of the flow in
the gaps between the boundary and the nearby inclusions. Roughly
speaking, the mean of the normal flow at the boundary enters the
inclusions, and thus excites the resistor network. The remainder, the
oscillations about the mean, have no effect on the network, but they
may be strong, depending on $k$ and $\delta$. As we explain below, the
term $\cR_k$ is important only in a specific ``resonant'' regime.

We distinguish three asymptotic regimes based on the values of the
dimensionless parameters
\begin{equation}
\ep = \frac{k \de}{L}\,, \qquad \eta = \frac{k R}{L}\, . 
\label{eq:PARA1}
\end{equation}
Equation (\ref{eq:F13}) implies that 
\begin{equation}
\ep \ll \eta\, \,
\label{eq:PARA2}
\end{equation}
but depending on the value of $k$, these parameters may be large 
or small. 

In the first regime $k \lesssim L/R$, so that 
\begin{equation}
\varepsilon \ll \eta \lesssim 1\, .
\label{eq:Reg1.1}
\end{equation}
The network is excited in this regime, and equation (\ref{eq:RN3})
shows that its boundary potentials $\Psi_i$ are basically the point
values of $\psi$ at the boundary nodes $\bx_i^\Gamma$. The energy of
the network plays an important role in the approximation, and it is
large, given by the sum of terms proportional to the effective
conductivities $\os_i$ and $\os_{ij}$ of the gaps, which are $O\left(
  \sqrt{R/\de} \right)$.  The term $\cR_k$ is much smaller, as
obtained from the following asymptotic expansions of the exponential
\begin{equation}
e^{-\sqrt{\ep \eta}} = 1 - \sqrt{\ep \eta} + O(\ep \eta)\, ,
\end{equation}
and the Polylogarithm function
\begin{equation}
\label{eq:Poly} 
{\rm Li}_{1/2}\left(e^{-2 \ep }\right) = \sqrt{\frac{\pi}{2 \ep}} +
\zeta\left(\frac{1}{2}\right) - 2 \ep \zeta\left(-\frac{1}{2}\right)
+ O\left(\ep^{3/2}\right)\, , \quad \ep \ll 1\, ,
\end{equation}
where $\zeta$ is the Riemann zeta function. We obtain that
\begin{align}
\label{eq:Reg1.2}
\cR_{i,k} =  \os_i\, O\left(\ep^{1/2}\right) \ll
\os_i\, , 
\end{align}
and conclude that $\cR_k$ is negligible in this regime. The
leading order of the energy is given by
\begin{equation}
  \cE(\psi) \approx E^{\rm net} \left( 
\bpsi(\psi) \right) + \frac{k \pi}{2} \,.
\label{eq:Reg1.3}
\end{equation}

In the second regime the boundary potential is very oscillatory, 
with $k \gtrsim L/\de \gg 1$, so that 
\begin{equation}
1 \lesssim \ep \ll \eta.
\label{eq:Reg2.1}
\end{equation}
The network plays no role in this regime, because it is not excited.
Its boundary potentials are exponentially small, essentially zero, as
shown in equation (\ref{eq:RN3}). The term $\cR_k$ is the sum of
\begin{align}
  \cR_{i,k} &= \frac{\os_i}{4} \sqrt{\frac{2 k \de_i}{\pi L}} e^{-
    \frac{2 k \de_i}{L} } \left[ 1 + O\left(e^{-
        \frac{2 k \de_i}{L} } \right) \right] \nonumber \\
  &= \frac{1}{2}\sqrt{\frac{\pi k R}{2L}} e^{-2 \ep \de_i/\de } \left[ 1
    + O\left(e^{- 2 \ep } \right) \right] \, ,
\end{align}
where we used the asymptotic expansion of the Polylogarithm function
at small arguments.  We estimate it as
\begin{equation}
\cR_k \sim \sqrt{\frac{k L}{R}} e^{-\ep}\, ,
\label{eq:residual}
\end{equation}
because 
\begin{equation}
N^\Gamma \sim \frac{L}{R}\, ,
\label{eq:NGamma}
\end{equation}
with symbol $\sim$ denoting henceforth approximate, up to a multiplicative
constant of order one. Consequently, 
\begin{equation}
\frac{\cR_k}{E_o(\psi)} \sim \sqrt{\frac{L}{Rk}} e^{-\ep}\, ,
\end{equation}
and recalling the definition (\ref{eq:PARA1}) of $\ep$, we see that
$\cR_k$ becomes negligible as $k$ increases. The oscillatory flow is
confined near the boundary $\Gamma$ for large $k$, and it does not see
the high contrast inclusions. The energy is approximately equal to
that in the reference medium
\begin{equation}
\cE(\psi) \approx E_o(\psi) = \frac{k \pi}{2}.
\label{eq:Reg2.2}
\end{equation}

The third regime corresponds to intermediate Fourier frequencies 
satisfying 
\[
\frac{L}{R} \lesssim  k \ll \frac{L}{\de},
\]
so that 
\begin{equation}
\ep \ll 1 \lesssim \eta.
\label{eq:Reg3.1}
\end{equation}
We call it the resonant regime because $\cR_k$ plays an important role
in the approximation.  Equation (\ref{eq:RN3}) shows that the network
gets excited, with boundary potentials that are smaller than the point
values of $\psi$.  The term $\cR_k$ is estimated by
\begin{equation}
  \cR_k \sim \sum_{i=1}^{N^\Gamma} \os_i \left[ 1 + O(\ep^{1/2}) 
  \right] \sim \frac{L}{R} 
\sqrt{\frac{R}{\de}} = \frac{k}{\sqrt{\ep \eta}}\, ,
\label{eq:Reg3.2}
\end{equation}
where we used the expansion (\ref{eq:Poly}), and (\ref{eq:NGamma}).
All the terms in (\ref{eq:R3}) play a role in the approximation of the
energy, with $\cR_k$ of the same order as $E^{\rm net}$ when $\ep \eta
\ll 1$, and much larger for $\ep \eta \gg 1$.  The term $\cR_k$
dominates the reference energy $E_o(\psi)$ when $\ep \eta \ll 1$, but
it plays a lesser role as the frequency $k$ increases so that $\ep
\eta \gtrsim 1$.

\subsection{General boundary potentials}
\label{sect:GB}
Assuming a potential of the form (\ref{eq:FDec}), with $K$ Fourier
modes, we write
\begin{equation}
\psi(\theta) = \sum_{k=0}^K \psi_k(\theta)\, ,
\label{eq:GB1}
\end{equation}
with $\psi_k(\theta)$ oscillating at frequency $k$, 
\begin{equation}
\psi_k(\theta) = a_k^c \cos(k \theta) + a_k^s \sin (k \theta)\, .
\label{eq:GB2}
\end{equation}
The maximum frequency $K$ may be arbitrarily large.  We obtain the
following generalization of the result in Theorem \ref{thm.1}.

\begin{corollary}
\label{cor.1}
For a potential $\psi$ of the form (\ref{eq:GB1}) we have that 
\begin{align}
  \left<\psi,\Lambda \psi \right> = 2 E(\psi) = 2 \left[ E^{\rm
      net}\left(\bpsi(\psi)\right) + \frac{1}{2}\left< \psi, \Lambda_o
      \psi\right> + \cR(\psi) \right] \left[ 1 + o(1)\right]\, .
\label{eq:GB3}
\end{align}
The first term is due to the network with boundary potentials
\begin{equation}
  \Psi_i(\psi) = \sum_{k=0}^K \psi_k(\theta_i) 
e^{-\frac{k \sqrt{2 R \de_i}}{L}}\, .
\label{eq:GB4}
\end{equation}
The second term is the quadratic form of the DtN map $\Lambda_o$ of
the reference medium, with uniform conductivity $\sigma_o = 1$. The
last term $\cR$ is given by 
\begin{align}
  \cR = \sum_{i=1}^{N^\Gamma} \sum_{k,m=0}^K e^{-|k-m| \frac{\sqrt{2 R
        \de_i}}{L}} \cR_{i,k\wedge m}\left\{\left(a_k^c a_m^c + a_k^s
      a_m^s\right)
    \cos[(k-m)\theta_i] + \right. \nonumber \\
  \left. \left(a_k^s a_m^c - a_k^c a_m^s\right)
    \sin[(k-m)\theta_i]\right\}\, , \label{eq:genRk}
\end{align}
where $ k \wedge m = \min\{k,m\}\, , $ and $\cR_{i,k}$ is defined in
(\ref{eq:defRik}).
\end{corollary}

The proof of this corollary is very similar to that of Theorem
\ref{thm.1}, so we do not include it here.  It uses the dual
variational principles (\ref{eq:F1}) and (\ref{eq:F7}) to estimate the
energy $E(\psi)$ for potential (\ref{eq:GB1}). Actually, it suffices
to consider
\[
\psi(\theta) = \cos(k \theta) + \cos(m \theta)\, , \quad 
\psi(\theta) = \sin(k \theta) + \cos(m \theta)\, , \quad 
\psi(\theta) = \sin(k \theta) + \sin(m \theta)\,,
\]
for arbitrary $k,m = 1, \ldots, K$, because the energy is a quadratic
form in $\psi$. We refer to \cite[section 4.3]{Thesis} for details.

The expression (\ref{eq:GB3}) is similar to (\ref{eq:R2}), and the
discussion in the previous section applies to the contribution of each
Fourier mode of $\psi$. The resonance $\cR$ captures the energy of the
oscillatory flow in the gaps between the inclusions and the boundary
$\Gamma$.  Its expression is more complicated than (\ref{eq:R11}), but
only the terms that are less oscillatory have a large contribution in (\ref{eq:genRk}).
  We can see this
explicitly in the special case where all the gaps are identical
\[
\de_i = \de_1\,, \quad \cR_{i,k} = \cR_{1,k}\, , \quad 
\forall \, i = 1, \ldots, N^\Gamma,
\]
and the boundary points are equidistant.  Then (\ref{eq:genRk})
simplifies to
\begin{align*}
  \cR =& 
  N^\Gamma  \sum_{k=0}^K\cR_{1,k}  
\left[\left(a_k^c\right)^{2}+ 
  \left(a_k^s \right)^{2}\right]  + \nonumber \\  
& 2 N^\Gamma  \sum_{k=0}^K\cR_{1,k} \sum_{q \in \mathbb{Z}^+} 
e^{-|q| N^\Gamma 
\frac{ \sqrt{2 R \de_1}}{L}}
1_{[0,K]}(k + q N^\Gamma)\left[a_k^c a_{k + q N^\Gamma}^c+ 
  a_k^s a_{k + q N^\Gamma}^s\right]  \, ,
\end{align*}
 because 
\[
\sum_{i=1}^{N^\Gamma}\cos[(k-m)\theta_i] = N^\Gamma \de_{km} ~{\rm modulo }\, N^\Gamma\, , 
\quad \sum_{i=1}^{N^\Gamma}\sin[(k-m)\theta_i] = 0\, , 
\quad \theta_i = \frac{(i-1)2 \pi}{N^\Gamma}\, .
\]
Here we let $1_{[0,K]}$ be the indicator function of the
interval $[0,K]$.

\subsection{Generalization to inclusions of different size and shape}
\label{sect:GENRadii}
We assumed for simplicity of the analysis that the inclusions $\cD_i$
are identical disks of radius $R$, but the results extend easily to
inclusions of different radii and even shapes. The leading order of
the energy is due to the singular flow in the gaps between the
inclusions and near the boundary. As long as we can approximate the
boundaries $\partial \cD_i$ locally, in the gaps, by arcs of circles
of radius $R_i$, and we have the scale ordering
\[
\de \sim \de_i \ll R_i \sim R \ll L,
\] 
the results of Theorem \ref{thm.1} and Corollary \ref{cor.1} apply,
with the following modifications: The effective conductivities of the
gaps are given by
\begin{equation}
  \os_{ij} = \pi \sqrt{\frac{2 R_i R_j}{\de_{ij}(R_i+R_j)}}\, , \quad 
i = 1, \ldots, N\, , ~~ j \in \cN_i\, ,
  \label{eq:osij_m}
\end{equation}
 and 
\begin{equation}
  \os_{i} = \pi \sqrt{\frac{2 R_i}{\de_i}}\, , \quad i = 1, 
\ldots, N^\Gamma.
\end{equation}
The resonance terms have the same expression as in (\ref{eq:R11}) and
(\ref{eq:genRk}), but $R$ is replaced by the local radii $R_i$ of
curvature in the sum over the gaps.

\section{Method of proof}
\label{sect:proofs}
The basic idea of the proof is to use the two variational principles
(\ref{eq:F1}) and (\ref{eq:F7}), with carefully chosen test potentials
and fluxes, to obtain upper and lower bounds on the energy that match
to the leading order, uniformly in $k$. The main difficulty in the
construction of these test functions is that, depending on $k$, the
flow may have very different behavior near $\Gamma$ than in the
interior of the domain. To mitigate this difficulty, we borrow an idea
from \cite{BGN-09,Novikov-09} and introduce in section \ref{sect:Perf} an
auxiliary problem in a so-called perforated domain $\Omega_p$. It is
a subset of $\Omega$, with complement $\Omega \setminus \Omega_p$
chosen so that the flow in it is diffuse, and thus negligible to
leading order in the calculation of energy.

The perforated domain is the union of two disjoint sets: the boundary
layer $\cB$, and the union of the gaps between the inclusions, denoted
by $\Pi$.  It is useful because it allows us to separate the analysis
of the energy in the boundary layer and that in $\Pi$, as shown in
section \ref{sect:PerfAdv}.  The estimation of the energy in $\Pi$ is
in section \ref{sect:Rnet}, where we review the network approximation.
The energy in $\cB$ is estimated in section \ref{sect:BLay}. The proof
of Theorem \ref{thm.1} is finalized in section \ref{sect:Est}.

\subsection{The perforated domain}
\label{sect:Perf}
\begin{figure}[t]
\begin{center}
\includegraphics[width=0.5\textwidth]{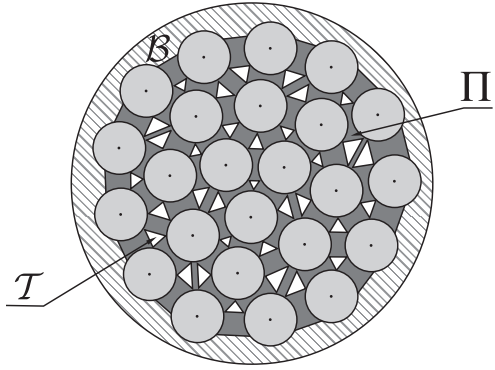}
\end{center}
\caption{Illustration of the perforated domain $\Omega_p$.  It
  is the union of two disjoint sets: the boundary layer $\cB$ and the
  set $\Pi$ of gaps between the adjacent inclusions.  The complement
  of $\Omega_{p}$ in $\Omega$ is the set $\cT$ of triangles.}
\label{fig:perf}
\end{figure}
\begin{figure}[ht]
\centering
\subfigure[]{   \includegraphics[scale=0.98]{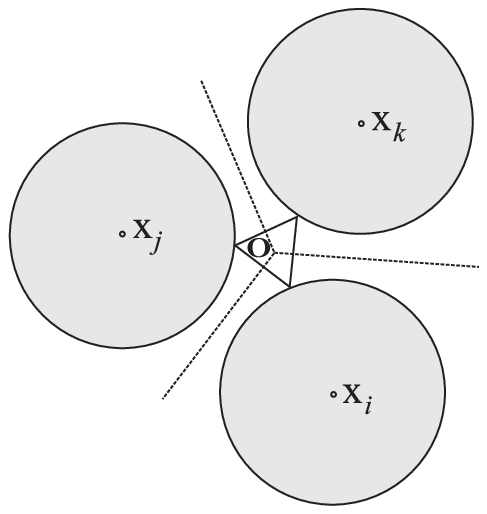} \label{fig:triangijk}  } \hspace{.3cm}
\subfigure[]{   \includegraphics[scale=1.0]{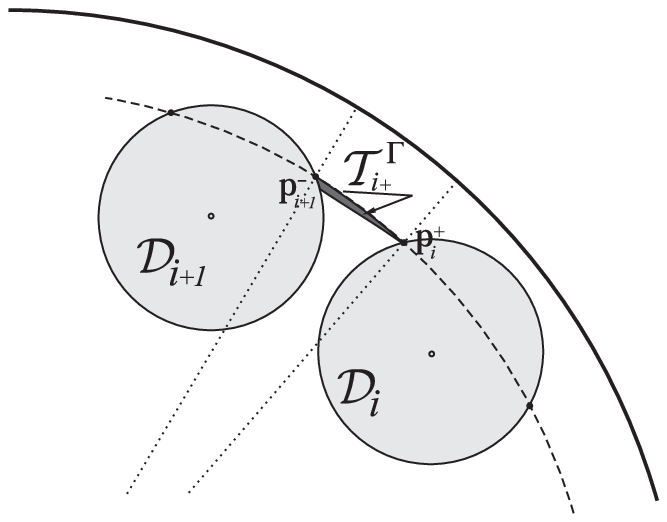} \label{fig:triangB} }
\label{F:domain-Liliana}
\caption{(a) Illustration of a triangle $\cT_{ijk}$. Its vertices are the
  intersections of the boundaries of the inclusions with the line
  segments between their centers and the vertex ${\bf O}$ of the
  Voronoi tessellation. \hspace{.015cm}  (b)  Illustration of a triangle $\cT_{i+}^\Gamma$. It has 
  vertices ${\bf p}_{i}^{+}$ and ${\bf p}_{i+1}^{-}$, one straight edge and two curved ones.
  One curved edge is the arc on the circle of radius $L-R/2$, shown with dashed line.
  The straight edge connects the vertex ${\bf p}_{i}^{+}$ with $\partial \cD_{i+1}$ along the 
  line that is parallel to that passing through $\bx_{i}$ and $\bx_{i+1}$. The other curved edge 
  is on $\partial \cD_{i+1}$.}
\end{figure}

Let us denote by $\cT$ the set of triangles that we wish to remove
from $\Omega$, based on the observation that the flow there is diffuse
and thus negligible in the calculation of the leading order of the
energy. There are two types of triangles, those in the interior of the
domain, and those near the boundary. The triangles in the interior are
denoted generically by $\cT_{ijk}$, for indexes $i \in \{1, \ldots,
N\}$, $j \in \cN_i$ and $k \in \cN_j$. We illustrate one of them in
Figure \ref{fig:triangijk}, where we denote by ${\bf O}$ the vertex of
the Voronoi tessellation, the intersection of the Voronoi cells
\[ {\bf O} = {\mathfrak{X}_i} \bigcap {\mathfrak{X}_i}\bigcap
{\mathfrak{X}_k}\, .
\]
The vertices of the triangle $\cT_{ijk}$ are at the intersections of
the boundaries of the inclusions with the line segments connecting
their centers with ${\bf O}$.

The triangles near the boundary are denoted by $\cT_{i+}^\Gamma$, for
$i = 1, \ldots, N^\Gamma$. Note that with our counting of the
inclusions  the triangle $\cT_{i+}^\Gamma$ involves the neighbors
$\cD_i$ and $\cD_{i+1}$ for $i = 1, \ldots, N^\Gamma-1$, whereas
$\cT^\Gamma_{N^\Gamma+}$ involves $\cD_{N^\Gamma}$ and $\cD_1$.  
We define the triangles to have one straight edge and two curved ones.
Let ${\bf p}_{i}^{\pm}$ be the intersection of the circle\footnote{
The circle of 
radius $L-R/2$ used in the definition of $\cT^\Gamma$ is somewhat
arbitrary.  We may chose any radius $L -R+\rho$, with $\sqrt{R \de}
\ll \rho \lesssim R/2$ and the result would be the same to the leading
order.} of radius $L-R/2$ shown with the dashed line in Figure \ref{fig:triangB} and the boundary $\partial \cD_{i}$ of the $i-$th inclusion. Then ${\bf p}_{i}^{+}$ and ${\bf p}_{i+1}^{-}$ are vertices
of $\cT_{i+}^{\Gamma}$ and the arc of the circle of radius $L-R/2$ between them is one curved 
edge of $\cT_{i+}^{\Gamma}$.  To determine the  straight edge of $\cT_{i+}^\Gamma$, we draw two  line segments that 
are parallel to the line through the centers $\bx_{i}$ and $\bx_{i+1}$ of the inclusions, and connect 
${\bf p}_{i}^{+}$ with $\partial \cD_{i+1}$ and ${\bf p}_{i+1}^{-}$ with $\partial \cD_{i}$, respectively. 
One of these segments lies inside the circle of radius $L-R/2$,  and it is the straight edge of $\cT_{i+}^{\Gamma}$.
The remaining curved edge is an arc on the boundary of one of the inclusions. If the straight 
edge stems from ${\bf p}_{i}^{+}$ i.e.,  if  $\cD_{i}$ is closer to $\Gamma$ 
than $\cD_{i+1}$, the curved edge lies on $\partial \cD_{i+1}$, as illustrated in Figure \ref{fig:triangB}. Otherwise 
it lies on $\partial \cD_{i}$. In the special case where the  two inclusions have the same distance to the boundary $\Gamma$, 
this  edge degenerates to a point. That is to say, $\cT_{i+}^ \Gamma $ has only two vertices 
${\bf p}_{i}^{+}$ and ${\bf p}_{i+1}^{-}$, and two edges connecting them. One edge is straight and the other is on the circle of radius $L-R/2$.  

The perforated domain is defined by 
\begin{equation}
\Omega_p = \Omega \setminus \cT.
\label{eq:perfD}
\end{equation}
It is the union of the boundary layer $\cB$ and the set of gaps $\Pi$, as shown 
in Figure \ref{fig:perf}..
The set $\cB$ is bounded on one side by $\Gamma$, and on the other
side by the inclusion boundaries $\partial \cD_i$ and the curved edges
of the triangles $\cT_{i+}^\Gamma$ between them, for $i = 1, \ldots,
N^\Gamma$. The set $\Pi$ is the union of the disjoint gaps $\Pi_{ij}$
between neighboring inclusions
\begin{equation}
\Pi = \bigcup_{
i = 1, \ldots, N, j \in \cN_i} \Pi_{ij}.
\label{eq:gaps}
\end{equation}
They are bounded by $\partial \cD_i$, $\partial \cD_j$, and the edges
of the interior triangles.

\subsection{Advantage of the perforated domain}
\label{sect:PerfAdv}

We define the energy $E_p(\psi)$ in the perforated domain by
\begin{equation}
E_p(\psi) = \min_{v \in \mathbb{V}_p(\psi)} \frac{1}{2} 
\int_{\Omega_p} d \bx \, |\nabla v(\bx)|^2\, ,
\label{eq:pVar}
\end{equation}
where the minimization is over potentials in the function space
\begin{equation}
  \mathbb{V}_p(\psi) = \left\{ v \in H^1(\Omega_p), ~ ~ v|_\Gamma 
    = \psi, ~~ v|_{\partial \cD_i} = ~{\rm constant}, ~ i = 1, 
    \ldots, N \right\}\, .
\label{eq:pVarF}
\end{equation}
Note that the set $\mathbb{V}(\psi)$ of test potentials in the
variational principle (\ref{eq:F1}) of $E(\psi)$ is contained in
$\mathbb{V}_p(\psi)$. Note also that the minimizer in (\ref{eq:pVar})
is the solution $u_p(\bx)$ of the Euler-Lagrange equations
\begin{align}
  \Delta u_p(\bx) &=0\, , \qquad ~\bx \in \Omega_p\, , 
\label{eq:pE1} \\
  u_p(\bx) &= \cU_i\,, \quad ~ ~\bx \in \partial \cD_i\, , 
\label{eq:pE2} \\
  \int_{\partial \cD_i} ds(\bx) \, {\bf n}(\bx) \cdot \nabla u_p(\bx)
  &=0\,, \qquad ~ i = 1, \ldots, N\,, \label{eq:pE3}\\
  u_p(\bx) &= \psi(\bx)\, , ~ ~ \bx \in \Gamma \,,
  \label{eq:pE4} \\
  {\bf n}(\bx) \cdot \nabla u_p(\bx) &=0\,, \qquad ~ \bx \in \partial
  \cT\, . \label{eq:pE5}
\end{align}
The first four equations are the same as those satisfied by the
minimizer of (\ref{eq:F1}), except that $\Omega_p$ is a subset of
$\Omega$. The unknowns are $u_p$ and the vector $\bU = (\cU_1, \ldots,
\cU_N)$ of constant potentials on the inclusions, the Lagrange
multipliers for the conservation of currents conditions
(\ref{eq:pE3}). Equation (\ref{eq:pE5}) says that there is no flow in
the set $\cT$ of triangles removed from $\Omega$.  The minimizer
$u(\bx)$ in (\ref{eq:F1}) does not satisfy these conditions, so
\[
u_p(\bx) \ne u(\bx), \quad \bx \in \Omega_p.
\]
However, the next lemma states that when replacing $u$ with $u_p$ we
make a negligible error in the calculation of the energy. The proof is 
in appendix \ref{ap:proofLem1}.

\begin{lemma}
\label{lem.1}
The energy $E(\psi)$ is approximated to leading order by the energy in
the perforated domain, uniformly in $k$,
\begin{equation}
E(\psi) = E_p(\psi) \left[ 1 + o(1)\right]\, .
\label{eq:lem1}
\end{equation}
\end{lemma}

Because the perforated domain is the union of the disjoint sets $\cB$
and $\Pi$, it allows us to separate the estimation of the energy in
the boundary layer from that in the gaps, as stated in the next lemma.
The two problems are tied together by the vector $\bU^\Gamma = (\cU_1,
\ldots, \cU_{N^\Gamma})$ of potentials on the inclusions near
$\Gamma$.

\begin{lemma}
\label{lem.2}
The energy in the perforated domain is given by the iterative
minimization
\begin{equation}
E_p(\psi) = \min_{\bU^\Gamma \in \mathbb{R}^{N^\Gamma}} 
\left[ E_\cB(\bU^\Gamma,\psi) + E_\Pi(\bU^\Gamma)\right]\, ,
\label{eq:lem2}
\end{equation}
where $E_\cB(\bU^\Gamma,\psi)$ and $E_\Pi(\bU^\Gamma)$ are the energy
in the boundary layer and gaps respectively, for given $\bU^\Gamma$
and $\psi$. The energy in the boundary layer has the variational
principle
\begin{equation}
  E_\cB(\bU^\Gamma,\psi) = \min_{v \in \mathbb{V}_\cB(\bU^\Gamma,\psi)}
\frac{1}{2} \int_\cB d \bx \, |\nabla v(\bx)|^2\, , 
\label{eq:EB}
\end{equation}
with minimization over potentials in the function space
\begin{align}
  \mathbb{V}_\cB(\bU^\Gamma,\psi) = \left\{ v \in H^1(\cB), ~ ~
    v|_\Gamma = \psi, ~ ~ v|_{\partial \cD_i} =
    \cU_i, \quad i = 1, \ldots, N^{\Gamma}, \right. \nonumber \\
  \left.  v|_{\partial \cD_i} = ~{\rm constant}, ~ i = N^{\Gamma}+1,
    \ldots, N \right\} \, .
\label{eq:VEB}
\end{align}
The energy in the gaps is given by 
\begin{equation}
E_\Pi(\bU^\Gamma) = \min_{v \in \mathbb{V}_\Pi(\bU^\Gamma)}
\frac{1}{2} \int_\Pi d \bx \, |\nabla v(\bx)|^2\, , 
\label{eq:EPi}
\end{equation}
with potentials in the function space
\begin{align}
  \mathbb{V}_\Pi(\bU^\Gamma) = \left\{ v \in H^1(\Pi), ~ ~
    v|_{\partial \cD_i} =
    \cU_i, \quad i = 1, \ldots, N^{\Gamma}, \right. \nonumber \\
  \left.  v|_{\partial \cD_i} = ~{\rm constant}, ~ i = N^{\Gamma}+1,
    \ldots, N \right\} \, .
  \label{eq:VEPi}
\end{align}
\end{lemma}

The proof of this lemma is in appendix \ref{ap:proofLem2}. It uses
that the minimizer $u_\cB$ of (\ref{eq:EB}) satisfies the
Euler-Lagrange equations
\begin{align}
  \Delta u_\cB(\bx) &= 0\,, \qquad ~ \bx \in \cB\, , \label{eq:B1} \\
  u_\cB(\bx) &= \cU_i\, , \quad ~ ~~ \bx \in \cD_i \, , ~ i = 1, \ldots,
  N^\Gamma\, , \label{eq:B2} \\
  u_\cB(\bx) &= \psi(\bx)\, , ~~ \bx \in \Gamma\, ,
  \label{eq:B3} \\
  {\bf n}(\bx) \cdot \nabla u_\cB(\bx) &=0\, , \qquad ~ 
\bx \in \partial \cB \cap \partial \cT\, , \label{eq:B4} 
\end{align}
and the minimizer $u_\Pi$ of (\ref{eq:EPi}) satisfies
\begin{align}
  \Delta u_\Pi(\bx) &= 0\,, \qquad ~ \bx \in \Pi\, , \label{eq:Pi1} \\
  u_\Pi(\bx) &= \cU_i\, , \quad ~ ~ \bx \in \cD_i \, , ~i = 1, \ldots,
  N\, , \label{eq:Pi2} \\
  \int_{\partial \cD_i} d s(\bx) \, {\bf n}(\bx) \cdot \nabla u_\Pi(\bx)
  &=0\, , \qquad ~ i = N^{\Gamma}+1, \ldots, N\,, \label{eq:Pi3} \\
  {\bf n}(\bx) \cdot \nabla u_\Pi(\bx) &=0\, , \qquad ~ 
\bx \in \partial \Pi \cap \partial \cT\, . \label{eq:Pi4} 
\end{align}
These equations are similar to (\ref{eq:pE1})-(\ref{eq:pE5}). Note
however that in (\ref{eq:B1})-(\ref{eq:B4}) there is only one unknown,
the potential function $u_\cB(\bx)$. The constant potentials on the
inclusions near the boundary are given. We do not get conservation of
current at the boundaries of these inclusions until we minimize
(\ref{eq:lem2}) over the vector $\bU^\Gamma$.  The unknowns in
equations (\ref{eq:Pi1})-(\ref{eq:Pi4}) are the potential function
$u_\Pi(\bx)$ and the vector $(\cU_{N^\Gamma+1}, \ldots, \cU_N)$ of
potentials on the interior inclusions. There is no explicit dependence
of $u_\Pi$ on the boundary potential $\psi$. The dependence comes
through $\bU^\Gamma$, when we minimize (\ref{eq:lem2}) over it.

We estimate in the next two sections $E_\Pi(\bU^\Gamma)$ and
$E_\cB(\bU^\Gamma,\psi)$. Then we gather the results and complete 
the proof of Theorem \ref{thm.1} in section \ref{sect:Est}.
\subsection{Energy in the gaps}
\label{sect:Rnet}
The energy $E_\Pi(\bU^\Gamma)$ is given by (\ref{eq:EPi}). We follow
\cite{BGN-09,Novikov-09} and rewrite it in simpler form using that
$\Pi$ is the union of the disjoint gaps $\Pi_{ij}$, for $i = 1,
\ldots, N$ and $j \in \cN_i$.
\begin{lemma}
\label{lem.3}
The energy $E_\Pi(\bU^\Gamma)$ is given by the discrete minimization
\begin{equation}
  E_\Pi(\bU^\Gamma) = \min_{\bU^{I} \in \mathbb{R}^{N-N^\Gamma+1}} 
  \frac{1}{2}\sum_{i=1}^N \sum_{j\in \cN_i} (\cU_i-\cU_j)^2 E_{ij}\,,
\label{eq:G1}
\end{equation}
where 
\[
\bU^I = (\cU_{N^\Gamma+1},\ldots, \cU_N)
\]
is the vector of potentials on the interior inclusions and 
$E_{ij}$ is the normalized energy in the gap $\Pi_{ij}$. It is given 
by the variational principle 
\begin{equation}
E_{ij} = \min_{v \in \mathbb{V}_{ij}} \frac{1}{2} \int_{\Pi_{ij}}
d \bx \, |\nabla v(\bx)|^2\, , \label{eq:G2}
\end{equation}
where the minimization is over the function space of 
potentials
\begin{equation}
  \mathbb{V}_{ij} = \left\{v \in H^1(\Pi_{ij}), 
    ~ ~ v|_{\partial \cD_i}= \frac{1}{2}, ~ ~
    v|_{\partial \cD_j}= -\frac{1}{2} \right\}\, .
\label{eq:G3}
\end{equation}
\end{lemma}

\begin{figure}[t]
\begin{center}
\includegraphics[width=0.41\textwidth]{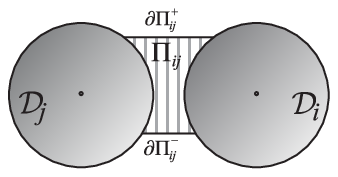}
\hspace{0.2in}
\raisebox{-0.1in}{\includegraphics[width=0.5\textwidth]{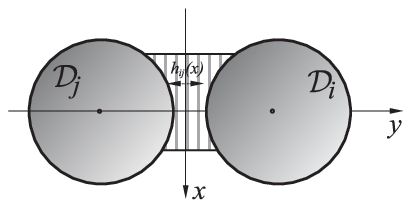}}
\end{center}
\caption{Illustration of a gap $\Pi_{ij}$. The local asymptotic analysis
is in the system of coordinates shown on the right, with $y$ axis 
connecting the centers of the inclusions.}
\label{fig:Pi}
\end{figure}

The proof of this iterative minimization is similar to that in
Appendix \ref{ap:proofLem2} and is given in \cite{BGN-09,Novikov-09}.
The estimate of the normalized energy $E_{ij}$ is obtained in
\cite{Berlyand-01,Keller-63}. It uses the variational principle
(\ref{eq:G2}) and a test potential $v(\bx)$ obtained from the
asymptotic approximation of the minimizer $u_{ij}(\bx)$ in the limit
$\de \to 0$ to obtain an upper bound of $E_{ij}$. The lower bound is
obtained from the dual variational principle
\begin{align}
  E_{ij} = \max_{{\bf j} \in \mathbb{J}_{ij}} \left[ \int_{\partial
      \cD_i \cap \partial \Pi_{ij}} \hspace{-0.1in} d s(\bx)\,
    \frac{1}{2} {\bf n}(\bx) \cdot {\bf j}(\bx) + \int_{\partial \cD_j
      \cap \partial \Pi_{ij}} \hspace{-0.1in} d s(\bx)\,
    \left(-\frac{1}{2}\right)
    {\bf n}(\bx) \cdot {\bf j}(\bx) - \right. \nonumber \\
  \left.  \frac{1}{2} \int_{\Pi_{ij}} d \bx \, |{\bf
      j}(\bx)|^2\right]\, ,
\label{eq:G4}
\end{align}
with fluxes ${\bf j}$ in the function space
\begin{equation}
  \mathbb{J}_{ij} = \left\{ {\bf j} \in L^2(\Pi_{ij}), ~ ~ 
    \nabla \cdot {\bf j}=0~ ~ {\rm in} ~ \Pi_{ij}, ~ ~ 
    {\bf n} \cdot {\bf j} = 0 ~ ~ {\rm in} ~ ~ 
    \partial \Pi_{ij}^\pm \right\} \, .
\label{eq:G5}
\end{equation}
Here $\partial \Pi_{ij}^\pm$ are the boundaries shared by $\Pi_{ij}$ and 
the set $\cT$ of triangles, as shown on the left in Figure 
\ref{fig:Pi}. 

The minimizing potential $u_{ij}$ of (\ref{eq:G2}) satisfies 
\begin{align}
\Delta u_{ij}(\bx) &=0\, , \quad ~~ \bx \in \Pi_{ij}\, , \\
u_{ij}(\bx) &=\frac{1}{2}\, , \quad ~~ \bx \in \partial \cD_i\, , \\
u_{ij}(\bx) &=-\frac{1}{2}\, , ~ ~ ~  \bx \in \partial \cD_j\, , \\
{\bf n}(\bx) \cdot \nabla u_{ij}(\bx) &=0\,, \quad ~~~ \bx \in \partial 
\Pi_{ij}^\pm\, .
\end{align} 
In the system of coordinates shown in Figure \ref{fig:Pi}, with $\bx =
(x,y)$ and $y $ axis along the line connecting the centers of the
inclusions, we see that $x$ belongs to an interval of order $R$ and
$y$ belongs to an interval of length
\begin{equation}
  h_{ij}(x) = \de_{ij} + 2 R \left(1 - \sqrt{1-\frac{x^2}{R^2}}\right)
  \, .
\label{eq:hij} 
\end{equation}
We expect that the leading order contribution to the energy comes from
the center of the gap, where $h_{ij} \sim \de \ll R$ and the gradient
of the potential is high. A simple scaling argument shows
that we can approximate $u_{ij}$ there by the potential $v$ satisfying
\[
\partial_{y}^2 v(x,y) = 0,
\]
with boundary conditions $v(x,\pm h_{ij}(x)/2) = \pm 1/2$. We obtain  
as in \cite{Berlyand-01,Keller-63}
\begin{equation}
v(\bx) = \frac{y}{h_{ij}(x)} \, ,
\label{eq:G6}
\end{equation}
and let the test flux be the divergence free vector that is
approximately equal to its gradient
\begin{equation}
{\bf j}(\bx) = \frac{1}{ h_{ij}(x) }{\bf e}_y \, . 
\label{eq:G7}
\end{equation}
Here ${\bf e}_y$ is the unit vector along the $y$ axis. It is parallel
to the boundaries $\partial \Pi_{ij}^\pm$ by construction, so
(\ref{eq:G7}) satisfies the no flow conditions there.

It is shown in \cite{Berlyand-01,Keller-63} that the upper bound
obtained with the test potential (\ref{eq:G6}) is given by
\begin{align}
  \frac{1}{2}\int_{\Pi_{ij}} d \bx\, \left|\nabla v(\bx)\right|^2 &=
  \frac{1}{2} \int_{-R}^R dx
  \int_{-\frac{h_{ij}(x)}{2}}^{\frac{h_{ij}(x)}{2}} dy\,
  \frac{1}{h^2_{ij}(x)} + O(1) \nonumber \\
   & = \frac{1}{2} \int_{-R}^R
  \frac{dx}{h_{ij}(x)} + O(1) \nonumber \\
  &= \frac{\os_{ij}}{2} + O(1)\, ,
\label{eq:OLDINT}
\end{align}
with 
\[
\os_{ij} = {\pi} \sqrt{\frac{R}{\de_{ij}}}\, .
\]
Moreover, the difference between the upper bound and the lower bound
given by the test flux (\ref{eq:G7}) is order one. Therefore,
\begin{equation}
  E_{ij} =  
  \frac{\os_{ij}}{2} + O(1)\, \, 
\label{eq:G8}
\end{equation}
and the energy $E_\Pi(\bU^\Gamma)$ follows from Lemma \ref{lem.3}
\begin{equation}
E_\Pi(\bU^\Gamma) = \cE_\Pi(\bU^\Gamma)\left[ 1 + o(1)\right]\, , 
\label{eq:G9}
\end{equation}
with leading order $\cE_\Pi$ given by
\begin{equation}
  \cE_\Pi(\bU^\Gamma) = \min_{\bU^{I} \in \mathbb{R}^{N-N^\Gamma+1}} 
  \frac{1}{2}\sum_{i=1}^N \sum_{j\in \cN_i} \frac{\os_{ij}}{2}(
  \cU_i-\cU_j)^2\, .
\label{eq:G10}
\end{equation}
This is the energy of the network with nodes at the centers $\bx_i$ of
the inclusions, edges $e_{ij}$ and net conductivities $\os_{ij}$, for
$i = 1, \ldots, N$ and $j \in \cN_i$. It is not the same network as in
Theorem \ref{thm.1}, because it does not contain the boundary nodes
$\bx_i^\Gamma$, for $i = 1, \ldots, N^\Gamma$.  It also has an
arbitrary vector $\bU^\Gamma$ of boundary potentials.  The network in
Theorem \ref{thm.1} has a uniquely defined vector $\bU^\Gamma(\psi)$
of potentials on the inclusions near $\Gamma$, the minimizer of
(\ref{eq:lem2}).

\subsection{Boundary layer analysis}
\label{sect:BLay}
To estimate the energy $E_\cB(\bU^\Gamma,\psi)$ we bound it above
using the variational principle (\ref{eq:EB}), and below using the
dual variational principle
\begin{align}
E_\cB(\bU^\Gamma,\psi) = \max_{{\bf j} \in \mathbb{J}_\cB} 
\left[ \int_\Gamma d s(\bx) \, {\bf n}(\bx) \cdot {\bf j}(\bx) + 
\sum_{i=1}^{N^\Gamma} \cU_i \int_{\partial \cB \cap \partial \cD_i} 
d s(\bx) \, {\bf n}(\bx) \cdot {\bf j}(\bx) - \right. \nonumber \\
\left. \frac{1}{2}
\int_{\cB} d\bx \, |{\bf j}(\bx)|^2 \right]\, ,
\label{eq:BL1}
\end{align}
with fluxes ${\bf j}$ in the function space
\begin{equation}
  \mathbb{J}_\cB = \left\{ {\bf j} \in L^2(\cB), ~ ~ 
    \nabla \cdot {\bf j} = 0 ~ {\rm in } ~ \cB, ~ ~ {\bf n} \cdot 
    {\bf j} = 0 ~ {\rm on} ~ \partial \cB \cap \partial \cT\right\}\, ,
\label{eq:BL2}
\end{equation}
and ${\bf n}$ the outer normal at $\partial \cB$.

Let us calculate the difference between the bounds, to gain insight in
the choice of the test potentials and fluxes in the variational
principles. We denote it by
\begin{align}
  \mathcal{G}(v,{\bf j}) =& \overline{E}_\cB(\bU^\Gamma,\psi;v) - 
\underline{E}_\cB(\bU^\Gamma,\psi;{\bf j})\, ,
\end{align}
 where
\begin{align}
  \overline{E}_\cB(\bU^\Gamma,\psi;v) = \frac{1}{2} \int_\cB d \bx \,
  |\nabla v(\bx)|^2\, , 
\end{align}
for $v \in \mathbb{V}_\cB(\bU^\Gamma,\psi)$  and 
\begin{align}
  \underline{E}_\cB(\bU^\Gamma,\psi;{\bf j})=\int_\Gamma d s(\bx) \,
  \psi(\bx) {\bf n}(\bx) \cdot {\bf j}(\bx) + \sum_{i=1}^{N^\Gamma}
  \cU_i \int_{\partial \cB \cap \partial \cD_i} d s(\bx) \, {\bf
    n}(\bx) \cdot {\bf j}(\bx) - \nonumber \\\frac{1}{2} \int_{\cB}
  d\bx \, |{\bf j}(\bx)|^2\, ,
\end{align}
for ${\bf j} \in \mathbb{J}_\cB$.  Integration by parts gives
\begin{align*}
  \int_\cB d \bx \, \nabla v(\bx) \cdot {\bf j}(\bx) &= \int_\cB d\bx
  \, \nabla \cdot \left[ v(\bx) {\bf j}(\bx)\right] \\
  &= \int_{\partial \cB} ds(\bx) \, v(\bx) {\bf n}(\bx) \cdot {\bf
    j}(\bx)\, \\
  &=\int_\Gamma ds(\bx) \, \psi(\bx) {\bf n}(\bx) \cdot {\bf j}(\bx) 
+ \sum_{i=1}^{N^\Gamma} \cU_i
  \int_{\partial \cB \cap \partial \cD_i} d s(\bx) \, {\bf n}(\bx)
  \cdot {\bf j}(\bx)\, ,
\end{align*}
because of the constraint $\nabla \cdot {\bf j} = 0$ and the boundary
conditions of $v$.  Therefore
\begin{equation}
  \mathcal{G}(v,{\bf j}) = \frac{1}{2} \int_\cB d \bx \, 
  |\nabla v(\bx) - {\bf  j}(\bx)|^2 \, ,
\label{eq:DefGap}
\end{equation}
and to make it small, we seek fluxes $ {\bf j}(\bx) \approx \nabla
v(\bx) $ in $\mathbb{J}_\cB$, and potentials $v \in
\mathbb{V}_\cB(\bU^\Gamma,\psi)$ satisfying
\begin{equation}
\label{eq:Dev}
\Delta v(\bx) \approx \nabla \cdot {\bf j}(\bx) = 0\, .
\end{equation}

\begin{figure}[t]
\begin{center}
\includegraphics[width=0.45\textwidth]{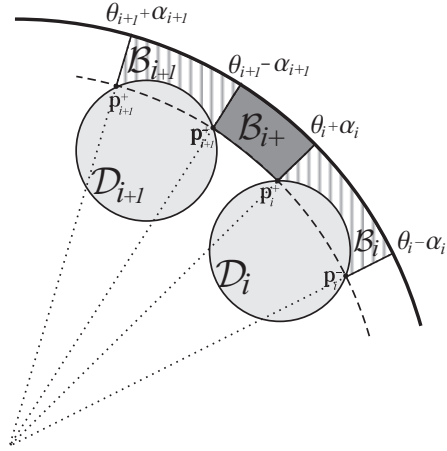}
\end{center}
\caption{Illustration of the decomposition of the boundary layer
  $\cB$.}
\label{fig:BLDec}
\end{figure}

\subsubsection{Test potentials for the upper bound}
\label{sect:TESTV}

Using the polar coordinates $(r,\theta)$ we write
\begin{equation}
  \cB = \{(r,\theta), ~ ~ r \in (L-d(\theta),L), ~ ~ 
\theta \in [0,2 \pi) \}\, ,
\end{equation}
with $d(\theta)$ the thickness of the layer given by 
\begin{equation}
  d(\theta) = \left\{ \begin{array}{ll} 
      \hspace{-0.07in}L-\rho_i \cos(\theta-\theta_i) - 
      \sqrt{R^2-\rho_i^2
        \sin^2(\theta-\theta_i)}\, ,  & \theta \in (\theta_i-
      \alpha_i,
      \theta_i+\alpha_i)\, , \\
      \hspace{-0.07in} \frac{R}{2}, & \theta \in 
      (\theta_i+\alpha_i,\theta_{i+1}-
      \alpha_{i+1})\,, ~ 
\end{array} \right.
\label{eq:dTh}
\end{equation}
where 
\begin{equation}
\rho_i = L - R - \de_i\, .
\end{equation}
The angles $\theta_{i} \pm \alpha_i$ are defined by the intersections
${\bf p}_i^\pm$ of the circle of radius $L-R/2$ with the boundaries
$\partial \cD_i$ of the inclusions. We estimate them as
\begin{equation}
  \sin \alpha_i \lesssim \frac{\sqrt{3}R}{2 \rho_i} = 
  O\left(\frac{R}{L}\right)\, , 
\label{eq:alpha}
\end{equation}
using Heron's formula for the triangle with edges of length $L-R/2$,
$\rho_i$ and $R$, and vertices at the origin, $\bx_i$ and ${\bf p}_i^+$.

Let us decompose the boundary layer in the sets
\begin{equation}
\cB_i = \{(r,\theta), ~ ~ r \in (L-d(\theta),L), ~ ~ \theta \in 
(\theta_i - \alpha_i,\theta_i+\alpha_i)\}\, , \quad i = 1, \ldots, 
N^\Gamma\, \,
\end{equation}
and 
\begin{equation}
  \cB_{i+} = \{(r,\theta), ~ ~ r \in (L-R/2,L), ~ ~ \theta \in 
  (\theta_i + \alpha_i,\theta_{i+1}-\alpha_{i+1})\}\, ,
\end{equation}
for $i = 1, \ldots, N^\Gamma-1$, as shown in Figure \ref{fig:BLDec}.
Recall that with our counting of the inclusions $\cD_1$ neighbors
$\cD_2$ and $\cD_{N^\Gamma}$, so we let
\begin{equation}
  \cB_{N^\Gamma+} = \{(r,\theta), ~ ~ r \in (L-R/2,L), ~ ~ \theta \in 
  (\theta_{N^\Gamma} + \alpha_{N^\Gamma},\theta_{1}-\alpha_{1})\} \,.
\end{equation}
We seek a test potential $v$ that is an approximate solution of
Laplace's equation in $\cB$, as stated in (\ref{eq:Dev}). We can solve
the equation with separation of variables in the domains $\cB_{i+}$,
but not in $\cB_i$, where the layer thickness varies with $\theta$.
However, the physics of the problem suggests that we neglect the
variation of $d(\theta)$ in the construction of $v$ in $\cB_i$.
Indeed, if it is the case that the tangential flow is dominant in
$\cB$, we expect that it is confined in a very thin layer near
$\Gamma$, of thickness smaller than $\de$, and does not interact with
the inclusions. Otherwise, the normal flow near $\Gamma$ plays a role,
and we expect that the leading contribution to the energy comes from
the gaps between the inclusions and $\Gamma$, where $d(\theta)$ is
smaller, of order $\de$.  Then, based on a scaling argument similar to
that in the previous section, we neglect the variation of $d(\theta)$
in the local approximation of the solution of Laplace's equation.

Consequently, we let the test potential be
\begin{align}
  v(r,\theta) = \left\{ \frac{(r/L)^k - [1-d(\theta)/L]^{2k} (L/r)^k}{
      1-[1-d(\theta)/L]^{2k}}\right\} \psi(\theta) +
  \frac{\ln(r/L)}{\ln[1-d(\theta)/L]}
  \mathfrak{L}(\theta,\bU^\Gamma)\, ,
\label{eq:testV}
\end{align}
where we recall that 
\[
\psi(\theta) = \cos (k \theta) \, .
\]
The function $\mathfrak{L}$ is constant on the inclusions
\begin{equation}
  \mathfrak{L}(\theta,\bU^\Gamma) = \cU_i\, , \quad 
\theta \in (\theta_i-\alpha_i,\theta_i+\alpha_i)\, ,
\end{equation}
and it interpolates linearly between the inclusions
\begin{equation}
\label{eq:interp}
\mathfrak{L}(\theta,\bU^\Gamma) = \frac{\cU_i+\cU_{i+1}}{2} 
+ \left(\cU_{i+1}-\cU_i\right) \left[\ell_i(\theta)-
  \frac{1}{2}\right]\, ,
\end{equation}
for $\theta \in (\theta_i+\alpha_i,\theta_{i+1}- \alpha_{i+1})$, where
\begin{equation}
  \ell_{i}(\theta) = \frac{\theta-(\theta_i + \alpha_i)}{(\theta_{i+1}-
    \alpha_{i+1})-(\theta_i+\alpha_i)}\, .
\end{equation}

The potential (\ref{eq:testV}) satisfies all the constraints in
$\mathbb{V}_\cB$, because
\begin{equation}
v|_\Gamma = v(L,\theta) = \psi(\theta)\, ,
\end{equation}
and 
\begin{equation}
  v|_{\partial \cD_i} =   v(L,d(\theta)) = 
  \mathfrak{L}(\theta,\bU^\Gamma) = 
  \cU_i, \quad \theta \in (\theta_i-\alpha_i, \theta_i+\alpha_i)\, .
\end{equation}
Thus, we can use it in the variational principle (\ref{eq:EB}) to
obtain an upper bound of the energy.

\subsubsection{Test fluxes for the lower bound}
\label{sect:TESTJ}
Since $v$ is harmonic by construction in the sets $\cB_{i+}$, we let 
\begin{equation} {\bf j}(r,\theta) = \nabla v(r,\theta) - {\bf e}_r
  \frac{(L-R/2)}{r}
  \partial_r v \left(L-\frac{R}{2},\theta\right) \quad \mbox{in} ~ 
\cB_{i+}\, ,
\label{eq:TJ1}
\end{equation}
where ${\bf e}_r$ is the unit vector in the radial direction.  We
obtain that
\[
\nabla \cdot {\bf j} = 0\quad \mbox{in} ~ \cB_{i+}\, ,
\]
and
\[ {\bf n} \cdot {\bf j}\left(L-\frac{R}{2},\theta \right) = -{\bf
  e}_r \cdot \nabla v \left(L-\frac{R}{2},\theta \right) + \partial_r
v \left(L-\frac{R}{2},\theta\right) = 0\,,
\]
for $\theta \in (\theta_i+\alpha_i,\theta_{i+1}-\alpha_{i+1})$, as
required by the constraints in $\mathbb{J}_\cB$.

Since in $\cB_i$ the potential $v$ is not harmonic, we cannot let the
flux be simply the gradient of $v$. We define it instead by
\begin{equation} {\bf j}(r,\theta) = \nabla^\perp H(r,\theta) =
  -\frac{{\bf e}_r}{r} \partial_\theta H(r,\theta) + {\bf e}_\theta
\partial_r H(r,\theta)\, ,
\label{eq:TJ1p}
\end{equation}
with scalar function 
\begin{equation}
H(r,\theta) = -\int_0^\theta d \theta' \, L \partial_r v(L,\theta') - 
\int_r^L \frac{dr'}{r'} \partial_\theta v(r',\theta)\, .
\label{eq:TJ2}
\end{equation}
This construction gives 
\[\nabla \cdot {\bf j}(\bx) = 0 \quad \mbox{in} ~ \cB_i\, ,
\]
with tangential flux equal to the tangential gradient of $v$ in
$\cB_i$
\[
  {\bf e}_\theta \cdot {\bf j}(r,\theta) = {\bf e}_\theta \cdot
  \nabla v(r,\theta)\, ,
\]
and normal flux matching the normal derivative of $v$ at $\Gamma$
\[
  {\bf e}_r \cdot {\bf j}(L,\theta) ={\bf e}_r \cdot \nabla
  v(L,\theta)\, .
\]
\subsubsection{The energy estimate}
\label{sect:Est}
We show in appendix \ref{ap:gap} that the test potential
(\ref{eq:testV}) and flux defined by (\ref{eq:TJ1}) and
(\ref{eq:TJ1p}) give the following difference between the upper and
lower bounds of $E_\cB$, 
\begin{equation}
\label{eq:GE1}
\mathcal{G}(v,{\bf j}) = \sum_{i =
  1}^{N^\Gamma} \left[\mathcal{G}_{\cB_i}(v,{\bf j}) + 
\mathcal{G}_{\cB_{i+}}(v,{\bf j})\right]\, ,
\end{equation}
where 
\begin{equation}
  \mathcal{G}_{\cB_i}(v,{\bf j})+
\mathcal{G}_{\cB_{i+}}(v,{\bf j}) \lesssim  O(1)\, .
\label{eq:GE2}
\end{equation}
The upper bound $\overline{E}_\cB(\bU^\Gamma,\psi;v)$ on the energy is
computed in appendix \ref{ap:BL_EN}.  We write it as
\begin{equation}
\label{eq:UB1}
\overline{E}_\cB(\bU^\Gamma,\psi;v) = \sum_{i =
  1}^{N^\Gamma} \left[\overline{E}_{\cB_i}(\bU^\Gamma,\psi;v) + 
\overline{E}_{\cB_{i+}}(\bU^\Gamma,\psi;v)\right]\, ,
\end{equation}
with terms
\begin{align}
  \overline{E}_{\cB_i}(\bU^\Gamma,\psi;v) &= \frac{1}{2} \int_{\cB_i}
  d\bx \, \left|\nabla v(\bx)\right|^2 \nonumber \\
  &= \frac{k \alpha_i}{2} + \frac{\pi}{2} \sqrt{\frac{2 L R}{\rho_i \de_i}}
  \left[\cU_i -\cos(k \theta_i) e^{-k \sqrt{\frac{2 R \de_i}{L
          \rho_i}}}\right]^2 + \nonumber \\& ~ ~ ~\frac{\pi}{4}
  \sqrt{\frac{2 L R}{\rho_i \de_i}} \left[ \sqrt{\frac{2 k \de_i}{L
        \pi}} {\rm Li}_{1/2} \left(e^{-2 k \de_i/L} \right)- e^{-k
      \sqrt{\frac{2 R \de_i}{L \rho_i}}}\right]+ O(1)\, ,
\label{eq:UB2} 
\end{align}
and 
\begin{align}
  \overline{E}_{\cB_{i+}}(\bU^\Gamma,\psi;v) &=\frac{1}{2}
  \int_{\cB_{i+}} d\bx \, \left|\nabla v(\bx)\right|^2 \nonumber \\& =
  \frac{k \left[(\theta_{i+1}-
      \alpha_{i+1})-(\theta_i+\alpha_i)\right]}{4} + O(1)\, .
\label{eq:UB3}
\end{align}
Note that the remainder is of the same order one as the difference
(\ref{eq:GE2}) between the upper and the lower bounds. We show next
that the first terms in (\ref{eq:UB2})-(\ref{eq:UB3}) are larger, and
thus define the leading order of the energy in $\cB$.

The magnitude of (\ref{eq:UB2})-(\ref{eq:UB3}) depends on the
potentials $\cU_i$ and the dimensionless parameters $\ep$ and $\eta$
defined in (\ref{eq:PARA1}), satisfying
\begin{equation}
  \ep = \frac{k \de}{L} \sim \frac{k \de_i}{L}\, , \qquad 
  \eta = \frac{k R}{L} \sim k \alpha_i \sim k (\theta_{i+1}-\theta_i)\, .
\end{equation} 
The potentials $\cU_i$ are arbitrary in (\ref{eq:UB2}), but in the end
we take them as minimizers of the energy, like in Lemma \ref{lem.2}.
They are the solutions of Kirchhoff's current conservation laws in the
network with boundary potentials $\cos(k \theta_i) e^{-k \sqrt{\frac{2
      R \de_i}{L \rho_i}}}$, for $i = 1, \ldots, N^\Gamma\, , $ and
satisfy the discrete maximum principle (\ref{eq:discMP}). Thus, we can
assume that
\[
|\cU_i -\cos(k \theta_i) e^{-k \sqrt{\frac{2 R \de_i}{L \rho_i}}}|
\sim e^{-k \sqrt{\frac{2 R \de_i}{L \rho_i}}} \sim e^{-\sqrt{\ep
    \eta}}\, .
\]

To see that the remainder of order one is negligible in (\ref{eq:UB2})
and (\ref{eq:UB3}), we distinguish two cases based on the value of
$\eta$. When $\eta \lesssim 1$, which means that $\ep \ll \eta
\lesssim 1$, the boundary potentials are of order one, the second term
in (\ref{eq:UB2}) dominates the others
\[
\frac{\pi}{2} \sqrt{\frac{2 L R}{\rho_i \de_i}} \left[\cU_i -\cos(k
  \theta_i) e^{-k \sqrt{\frac{2 R \de_i}{L \rho_i}}}\right]^2 \sim
\sqrt{\frac{R}{\de}} \gg O(1) \gtrsim k \alpha_i \sim \eta\,,
\]
and the remainder is negligible.  Otherwise, $\eta \gg 1$ and the
remainder is again negligible to leading order, because 
\[
k \alpha_i \sim  \eta \gg 1.
\]
We gather the results and rewrite the energy in the boundary layer as
\begin{align}
  E_\cB(\bU^\Gamma,\psi) = &\left\{ \frac{k \pi}{2} +
    \sum_{i=1}^{N^\Gamma} \frac{\pi}{2} \sqrt{\frac{2 L R}{\rho_i
        \de_i}} \left[\cU_i -\cos(k \theta_i) e^{-k \sqrt{\frac{2 R
            \de_i}{L
            \rho_i}}}\right]^2 + \right. \nonumber \\
  &\left. \frac{\pi}{4} \sqrt{\frac{2 L R}{\rho_i \de_i}} \left[
      \sqrt{\frac{2 k \de_i}{L \pi}} {\rm Li}_{1/2} \left(e^{-2 k
          \de_i/L} \right)- e^{-k \sqrt{\frac{2 R \de_i}{L
            \rho_i}}}\right]\right\} \left[1 + o(1)\right] \, ,
\label{eq:EB_R}
\end{align}
with negligible relative error, uniformly in $k$. 

The proof of Theorem \ref{thm.1} follows from (\ref{eq:EB_R}), once we
replace $\rho_i = L\left[1+O(R/L)\right]$ by $L$ in (\ref{eq:EB_R}).
We can do so without affecting the leading order, independent of the 
value of $k$. 

\section{Summary}
\label{sect:sum}
We obtained an asymptotic approximation of the Dirichlet to Neumann
(DtN) map $\Lambda$ of the partial differential equation describing
two dimensional electrical flow in a high contrast composite medium
occupying a bounded, simply connected domain $\cD$ with smooth
boundary $\Gamma$. The high contrast composite has perfectly
conducting inclusions packed close together, so they are close to
touching. To simplify the proofs, we assumed that $\cD$ is a disk of
radius $L$, and that the inclusions are identical disks of radius $R$.
Extensions to general domains, sizes and shapes of inclusions are
discussed, as well.  The analysis is in the regime of separation of
scales $\de \ll R \ll L$, where $\de$ is the typical thickness of the
gaps between adjacent inclusions.

Because the map $\Lambda$ is self-adjoint, it is determined by its
quadratic forms $\left<\psi,\Lambda \psi\right>$, for all boundary
potentials $\psi$ in the trace space $H^{1/2}(\Gamma)$. The main
result of the paper is the explicit characterization of the leading
order of these quadratic forms in the regime of separation of scales
described above.  The result is intuitive once we decompose the
potential $\psi$ over Fourier modes, and study the quadratic forms
$\left<\psi_k,\Lambda \psi_k\right>$ for modes $\psi_k$ oscillating at
arbitrary frequency $k$. It says that the leading order of
$\left<\psi_k,\Lambda \psi_k\right>$ is given by the sum of three
terms: The first is the quadratic form $\bpsi(\psi_k) \cdot
\Lambda^{\rm net} \bpsi(\psi_k)$ of the matrix valued DtN map
$\Lambda^{\rm net}$ of a unique resistor network with vector
$\bpsi(\psi_k)$ of boundary potentials. The second term is the
quadratic form $\left<\psi_k,\Lambda_o \psi_k\right>$ of the DtN map
$\Lambda_o$ of the homogeneous medium with reference conductivity
$\sigma_o = 1$ in which the inclusions are embedded. The last term
$\cR$ is labeled a resonance term, because it plays a role only in a
certain ``resonant'' regime.

The resistor network approximation arises due to the singularity of
the potential gradient in the gaps between the inclusions, and the
gaps between the boundary and the nearby inclusions. The network is
unique, with nodes at the centers of the inclusions and edges
connecting adjacent inclusions. The edge conductivities capture the
net energy in the associated gaps. Network approximations have been
derived before in homogenization studies of high contrast composites.
What is new here is that the excitation of the network, the vector of
potentials $\bpsi(\psi_k)$ at its boundary nodes, depends on the
frequency $k$ of oscillation of $\psi_k$. If $k$ is small, then the
entries in $\bpsi(\psi_k)$ are the values of $\psi_k$ at the points on
$\Gamma$ that are closest to the inclusions. However, for large $k$,
the entries in $\bpsi(\psi_k)$ are damped exponentially in $k$.  There
is a layer of strong flow near the boundary $\Gamma$, and the network
plays a lesser role as $k$ increases.  We distinguished three regimes
in the approximation of $\left<\psi_k,\Lambda \psi_k\right>$. In the
first regime $k$ is small, so that the entries in $\bpsi(\psi_k)$ are
large, of order one. The network is excited and plays a dominant role
in the approximation,
\[
\left<\psi_k,\Lambda \psi_k\right> \approx \bpsi(\psi_k) \cdot
\Lambda^{\rm net} \bpsi(\psi_k)\, .
\]
In the second regime the frequency $k$ is very large, so that the flow
is confined in a very thin layer near the boundary $\Gamma$ and does
not interact with the inclusions. The entries in $\bpsi(\psi_k)$ are
basically zero, the network is not excited and the flow perceives the
medium as homogeneous
\[
\left<\psi_k,\Lambda \psi_k\right> \approx 
\left<\psi_k,\Lambda_o \psi_k\right>\, .
\]
In the third, intermediary regime, some of the flow penetrates in the
domain and excites the network. The remainder is tangential flow near
the boundary, as in the homogeneous medium, and oscillatory flow
squeezed between the boundary and the nearby inclusions. The latter
gives an anomalous energy, captured by the resonance term $\cR$.  All
three terms play a role in the approximation in this resonant regime,
\[
\left<\psi_k,\Lambda \psi_k\right> \approx \bpsi(\psi_k) \cdot
\Lambda^{\rm net} \bpsi(\psi_k) + \left<\psi_k,\Lambda_o \psi_k\right>
+ 2\cR\, .
\]
Our analysis justifies these approximations and gives explicit
formulas for $\bpsi(\psi_k)$ and the resonant term $\cR$.

\section*{Acknowledgements}
The work of L. Borcea was partially supported by the AFSOR Grant
FA9550-12-1-0117, the ONR Grant N00014-12-1-0256 and by the NSF Grants
DMS-0907746, DMS-0934594.  The work of Y. Wang was supported by the
NSF Grants DMS-0907746, DMS-0934594. 
Y. Gorb was supported by the NSF grant DMS-1016531.

\appendix
\section{Maximum principle for the potentials on the inclusions}
\label{sect:MP}
We show here that the potentials $\cU_i$ on the inclusions $\cD_i$,
for $i = 1, \ldots, N$, are bounded in terms of the boundary data
$\psi$. 

Consider the solution $(u,\bU)$ of equations (\ref{eq:F3})-(\ref{eq:F6}). Since $u$ is harmonic in the connected
set $\Omega$, it takes its  maximum and minimum values at the boundary
$
\partial \Omega = \Gamma \displaystyle \bigcup_{i=1}^N \partial \cD_i\, .
$ 
Suppose that there exists an index $i$ for which 
\[
u|_{\partial \cD_{i}}= \cU_i = M = \max_{\bx \in \overline{\Omega}} u(\bx)\, ,
\]
and define the function
\begin{align}
  f(\rho) &= \frac{1}{2 \pi (R+\rho)} \int_{|\bx-\bx_i| = R + \rho}
  ds(\bx) \, u(\bx) \nonumber \\
  &=\frac{1}{2 \pi} \int_{|\by| = 1} ds (\by)\, u(\bx_i + (R+\rho)
  \by)\, , \label{eq:MP.1}
\end{align}
for $\rho \le O(\de)$, so that  the annulus
\[
\mathcal{C}_\rho = \left\{ \bx = \bx_i + r (\cos \theta,\sin 
\theta)\, , ~r \in [R,R+\rho]\, , ~ \theta \in [0,2\pi]\right\}\, ,
\]
is contained in $\Omega$.  We obtain using integration by parts 
and the conservation of currents (\ref{eq:F5}) at $\partial \cD_i$ that 
\begin{align*}
  f'(\rho) &= \frac{1}{2 \pi} \int_{|\by| = 1} ds (\by)\, \by \cdot
  \nabla u(\bx_i + (R+\rho) \by) \\
  &= \frac{1}{2 \pi(R + \rho)} \int_{|\bx-\bx_i| = R+\rho} 
ds (\bx) \, {\bf n(\bx)} \cdot \nabla u(\bx)
  \\
  &= \frac{1}{2 \pi(R + \rho)}\left[ \int_{|\bx-\bx_i| = R+\rho} ds
    (\bx)\, {\bf n(\bx)} \cdot \nabla u(\bx) + \int_{\cD_i}ds (\bx)\,
    {\bf n(\bx)} \cdot \nabla u(\bx) \right] \\
  &= \frac{1}{2 \pi(R + \rho)} \int_{\mathcal{C}_\rho} d\bx \,
  \Delta u(\bx) \\
  & = 0\, ,
\end{align*}
and therefore $f(\rho)$ is constant 
\begin{equation}
f(\rho) = f(0) = \cU_i = M\,.
\end{equation}
Moreover, integrating in polar coordinates we get that the average of 
$u(\bx)$ in the annulus equals its maximum value
\begin{align*}
  \frac{1}{|\mathcal{C}_\rho|}\int_{\mathcal{C}_\rho} d\bx \, u(\bx)
  &= \frac{1}{|\mathcal{C}_\rho|} \int_R^{R+\rho} dr \, 
\int_{|\bx-\bx_i| = r} ds(\bx) \, u(\bx) \\
&= \frac{1}{|\mathcal{C}_\rho|} \int_R^{R+\rho} dr 2 \pi r f(r-R) \\
&= M\, .
\end{align*} 
This implies that $ u(\bx) = M$ in $ \mathcal{C}_\rho$, and using the
maximum principle for the harmonic function $u(\bx)$, that $ u(\bx) =
M $, in $\Omega$.

A similar argument shows that if the minimum value of the potential is
attained at the boundary of one inclusion, then $u$ is constant in $\Omega$.
Thus, we have the maximum principle
\begin{equation}
\label{eq:MPR}
\min_{\bx \in \Gamma} \psi(\bx) \le \cU_{i}  \le \max_{\bx \in \Gamma} \psi(\bx), \quad 
i = 1, \ldots, N\, .
\end{equation}

A discrete version of the maximum principle for networks can be found
in \cite{RNet1,RNet2}. It says that the potential at the nodes of the
network attains its minimum and maximum values at the boundary nodes.
Thus, if we let $\Psi_i$ for $i = 1, \ldots, N^\Gamma$ be the boundary
potentials, we have 
\begin{equation}
  \min_{j = 1, \ldots,   N^\Gamma} 
  \Psi_j \le \cU_i \le \min_{j = 1, \ldots,   N^\Gamma} \Psi_j\, , 
  \quad  i = 1, \ldots, N^\Gamma\, .
\label{eq:discMP}
\end{equation}

\section{Proof of Lemma \ref{lem.1}}
\label{ap:proofLem1}
Recall that  the solution $u(\bx)$ of equations
(\ref{eq:F3})-(\ref{eq:F6}) minimizes (\ref{eq:F1}).  We have
\begin{equation}
E(\psi) = \frac{1}{2} \int_\Omega d \bx \, |\nabla u(\bx)|^2 \ge 
\frac{1}{2} \int_{\Omega_p} d \bx \, |\nabla u(\bx)|^2 \ge 
E_p(\psi)\, .
\label{eq:PL.1}
\end{equation}
The first inequality is because $\Omega_p \subset \Omega$. The second
inequality is because the restriction of $u(\bx)$ to $\Omega_p$
belongs to the function space $\mathbb{V}_p(\Psi)$ of test potentials
in the variational formulation (\ref{eq:pVar}) of $E_p(\Psi)$.  To
complete the proof we need the following result, obtained in sections
\ref{sect:Rnet} and \ref{sect:BLay}.

\begin{lemma}
\label{lem.A}
There exists a potential $v_p(\bx) \in \mathbb{V}_p(\psi)$ such that
\begin{equation}
\frac{1}{2}\int_{\Omega_p} d \bx \, |\nabla v_p(\bx)|^2 = E_p(\psi)
\left[ 1 + o(1)\right].
\label{eq:lemA}
\end{equation}
Moreover, if let $T$ be any edge of a triangle in $\cT$, and denote by
$|T|$ its length, we have the pointwise estimate
\begin{equation}
  |\nabla v_p(\bx) | \le \frac{C}{|T|}, \quad 
\bx \in T \subset \partial \cT\, , 
\label{eq:PL.2}
\end{equation}
with order one constant $C$ that is independent of $\de$ and $k$.
\end{lemma}

The estimate (\ref{eq:PL.2}) is valid in the vicinity of the
boundary of the triangles, not only on $\partial \cT$. Moreover, 
by our definition of the triangles,
\begin{equation}
|T| = O(R)\, .
\end{equation}
Using Kirszbraun's theorem \cite{federer-96} we extend $v_p(\bx)$ from
the boundary of each triangle inside the triangle, in such a way that
$|\nabla v_p|$ remains bounded by $O(1/R)$ in $\cT$.  The extended
$v_p$ is a function in $\mathbb{V}(\psi)$, so we get the upper bound
\begin{equation}
E(\psi) \le \frac{1}{2} \int_{\Omega} d\bx \, 
|\nabla v_p(\bx)|^2 = \frac{1}{2} \int_{\Omega_p} d\bx \, 
|\nabla v_p(\bx)|^2 + \frac{1}{2} \int_{\cT} d\bx \, 
|\nabla v_p(\bx)|^2\, .
\label{eq:PL.2p}
\end{equation}
The first term in the right hand side is given in (\ref{eq:lemA}).  To
estimate the second term, let $\cT_{ijk}$ be an arbitrary interior
triangle, for $i = 1,\ldots, N$, $j \in \cN_i$ and $k \in \cN_k$. By
construction, the area of the triangles is $O(R^2)$, so we have
\begin{equation}
\frac{1}{2} \int_{\cT_{ijk}} d\bx \, |\nabla v_p(\bx)|^2 = O(1)\, .
\label{eq:PL.3}
\end{equation}
This is much smaller than the  contribution of the
gaps given in section \ref{sect:Rnet}
\begin{equation}
  \frac{1}{2} \int_{\Pi_{ij}} d\bx \, |\nabla v_p(\bx)|^2 = 
O\left(\sqrt{\frac{R}{\de}}\right)\, , \quad 
  \frac{1}{2} \int_{\Pi_{jk}} d\bx |\nabla v_p(\bx)|^2 = 
O\left(\sqrt{\frac{R}{\de}}\right)\, .
\label{eq:PL.4}
\end{equation}
A similar result holds for the triangles near the boundary layer. 
We obtain that 
\begin{equation}
  \frac{1}{2} \int_{\cT} d\bx \, |\nabla v_p(\bx)|^2 = E_p(\psi) \, 
o(1)\, ,
\end{equation}
and the proof of Lemma \ref{lem.1} follows from (\ref{eq:PL.2p}) and
(\ref{eq:PL.1}).

\section{Proof of Lemma \ref{lem.2}}
\label{ap:proofLem2}
The proof is a consequence of Euler-Lagrange equations
(\ref{eq:pE1})-(\ref{eq:pE5}), (\ref{eq:B1})-(\ref{eq:B4}) and
(\ref{eq:Pi1})-(\ref{eq:Pi4}), which have unique solutions as follows
from standard application of Lax-Milgram's Theorem.  It is convenient
in this section to emphasize in the notation their dependence on the
data. Thus, we let $u_p(\bx;\psi)$, $\bU(\psi)$ be the solutions of
(\ref{eq:pE1})-(\ref{eq:pE5}). Moreover, for a given $\bU^\Gamma =
(\cU_1, \ldots, \cU_{N^\Gamma})$ we let $u_\cB(\bx;\psi,\bU^\Gamma)$
be the solution of (\ref{eq:B1})-(\ref{eq:B4}), and
$u_\Pi(\bx;\bU^\Gamma)$ and $\bU^{I}(\bU^\Gamma) =
(\cU_{N^\Gamma+1}(\bU^\Gamma), \ldots, \cU_{N}(\bU^\Gamma))$ the
solution of (\ref{eq:Pi1})-(\ref{eq:Pi4}). The index $I$ stands for
interior inclusions.

Note that the restriction of $u_p(\bx;\psi)$ to the boundary layer
solves equations (\ref{eq:B1})-(\ref{eq:B4}) for $\bU^\Gamma =
\bU^\Gamma(\psi) = ( \cU_1(\psi), \ldots, \cU_{N^\Gamma}(\psi))$,
\begin{equation}
u_\cB(\bx;\psi,\bU^\Gamma(\psi)) = u_p(\bx;\psi)\,,  \quad 
\bx \in \cB\,.
\label{eq:P2.1}
\end{equation}
Similarly, the restriction of $u_p(\bx;\psi)$ to the set $\Pi$ of 
gaps 
\begin{equation}
u_\Pi(\bx;\bU^\Gamma(\psi)) = u_p(\bx;\psi)\, , \quad 
\bx \in \Pi \, , 
\label{eq:P2.2}
\end{equation}
and the vector of potentials on the interior inclusions
\begin{equation}
  \bU^I(\bU^\Gamma(\psi)) = (\cU_{N^\Gamma+1}(\psi), 
\ldots, \cU_N(\psi))\, , \quad 
  \bx \in \Pi \, , 
\label{eq:P2.3}
\end{equation}
solve equations (\ref{eq:Pi1})-(\ref{eq:Pi4}) for $\bU^\Gamma =
\bU^\Gamma(\psi)$.
Therefore, we have 
\begin{align}
  E_p(\psi) &= \frac{1}{2} \int_{\Omega_p} d \bx \, | \nabla u_p (\bx
  ; \psi )|^2
  \nonumber \\
&= \frac{1}{2} \int_{\cB} d \bx \, | \nabla u_\cB (\bx; \psi ,
  \bU^\Gamma ( \psi )) |^2 + \frac{1}{2} \int_{\Pi} d \bx \, |\nabla
  u_\Pi ( \bx ;\bU^\Gamma ( \psi ))|^2 \nonumber \\
&=E_\cB( \bU^\Gamma ( \psi ),\psi ) +
  E_\Pi(\bU^\Gamma ( \psi )) \nonumber \\
&\ge \min_{\bU^\Gamma} \left[ E_\cB( \bU^\Gamma ,\psi ) + E_\Pi(
    \bU^\Gamma )\right]\, .
  \label{eq:P2.4}
\end{align}

For the reverse inequality let $\bU^\Gamma$ be arbitrary in
$\mathbb{R}^{N^\Gamma}$ and define $v \in \mathbb{V}_p(\psi)$ by 
\begin{equation}
  v(\bx) = \left\{ \begin{array}{ll}
u_\cB(\bx; \psi, \bU^\Gamma)\,,  &\quad \bx \in \cB\, , \\
u_\Pi(\bx;  \bU^\Gamma)\,,  &\quad \bx \in \Pi\, .
\end{array} \right.
\label{eq:P2.5}
\end{equation}
We obtain that 
\begin{align}
  E_p(\psi) &\le \frac{1}{2}\int_{\Omega_p} d \bx \,
  |\nabla v(\bx)|^2 \nonumber \\
  &= \frac{1}{2}\int_{\cB} d \bx \, |\nabla
  u_\cB(\bx;\psi,\bU^\Gamma)|^2 + \frac{1}{2}\int_{\Pi} d \bx \,
  |\nabla u_\Pi(\bx;\bU^\Gamma)|^2 \nonumber \\
&= E_\cB(\bU^\Gamma,\psi) + E_\Pi(\bU^\Gamma),
\label{eq:P2.6}
\end{align}
for all $\bU^\Gamma$. The result follows by taking the minimum
over $\bU^\Gamma$ in $\mathbb{R}^{N^\Gamma}$.

\section{Tightness of the bounds on $E_\cB$}
\label{ap:gap}
Definition (\ref{eq:TJ1}) of the flux in $\cB_{i+}$ and the expression
(\ref{eq:testV}) of the potential give that
\begin{align}
  \mathcal{G}_{\cB_{i+}}(v,{\bf j}) &= \frac{1}{2}\int_{\cB_{i+}} d
  \bx \, \left| \frac{(L-R/2)}{r} \partial_r v
    \left(L-\frac{R}{2},\theta\right)\right|^2 \nonumber \\
  &=\frac{1}{2} \int_{\cB_{i+}} d \bx \, \left|\frac{2 (L-R/2)}{rL}
    \frac{k \left(1 - \frac{R}{2 L}\right)^{k-1}}{1 - \left(1 -
        \frac{R}{2 L}\right)^{2k}}\cos (k \theta) +
    \frac{\mathfrak{L}(\theta,\bU^\Gamma)}{r \ln\left(1 - \frac{R}{2
          L}\right)}\right|^2 \, .
\label{eq:T1}
\end{align}
We estimate the first term by
\begin{align*}
\left| 
\frac{2 (L-R/2)}{rL} \frac{k \left(1 - \frac{R}{2 L}\right)^{k-1}}{
1 - \left(1 - \frac{R}{2L}\right)^{2k}}\cos (k \theta) \right| \le 
\frac{4}{R}\, , 
\end{align*}
because $r \in (L-R/2,L)$ and the function $k a^k/(1-a^{2k})$ for any
$a \in (0,1)$ is monotone decreasing in $k$ for $k \ge 1$. In
particular, for $a = 1 -R/(2L)$, we have
\[
\frac{k
      \left(1 - \frac{R}{2 L}\right)^{k}}{1 - \left(1 - \frac{R}{2
          L}\right)^{2k}} \le \frac{1- \frac{R}{2 L}}{1-1+\frac{R}{2 L}} = \frac{2L}{R} \left(1-\frac{R}{2L}\right)\, .
\]
The second term in (\ref{eq:T1}) satisfies 
\[
\left|\frac{\mathfrak{L}(\theta,\bU^\Gamma)}{r \ln\left(1 - \frac{R}{2
        L}\right)}\right| \le \frac{2}{R} \left(1 + O(R/L)\right)\, 
\]
because $\mathfrak{L}(\theta,\bU^\Gamma)$ is the interpolation between
$\cU_i$ and $\cU_{i+1}$, and their absolute value is bounded by one,
as shown in (\ref{eq:MPR}).  Thus, we have
\begin{equation}
\mathcal{G}_{\cB_{i+}}(v,{\bf j}) \le \frac{13
  \left(1+O(R/L)\right)}{R^2} \int_{\cB_{i+}} d \bx = O(1)\, , 
\end{equation}
because 
\begin{align*}
\int_{\cB_{i+}} d\bx = \int_{L-R/2}^{L} dr \, r \int_{\theta_{i}+
\alpha_{i}}^{\theta_{i+1}-\alpha_{i+1}} d\theta
= \frac{LR}{2} \left(1 - \frac{R}{4L}\right)\left(\theta_{i+1}-
\theta_{i}-\alpha_{i+1}-\alpha_{i}\right) \sim R^{2}\,,
\end{align*}
and 
\[
\theta_{i+1}-\theta_{i}-\alpha_{i+1}-\alpha_{i} \lesssim \frac{2
  \pi}{N^{\Gamma}} \sim \frac{R}{L}\, .
\]

Definition (\ref{eq:TJ1p}) of the test flux gives after a
straightforward calculation that
\begin{align}
  \mathcal{G}_{\cB_i}(v,\bj) &= \frac{1}{2} \int_{\cB_i} d \bx
  \left|\nabla v(\bx) - \bj(\bx)\right|^2 \nonumber \\
  &= \frac{1}{2} \int_{\theta_i-\alpha_i}^{\theta_i + \alpha_i} d
  \theta \int_{L-d(\theta)}^L dr \, r \left| \partial_r 
v(r,\theta) +
    \frac{1}{r} \partial_\theta H(r,\theta)\right|^2\, \nonumber \\
  &= \frac{1}{2} \int_{\theta_i-\alpha_i}^{\theta_i + \alpha_i} d
  \theta \int_{L-d(\theta)}^L \frac{dr}{r} \left[ \int_r^L ds \,s 
\Delta v(s,\theta)\right]^2\, ,
\end{align}
and using expression (\ref{eq:testV}) of the test potential 
we obtain 
\begin{align}
  \mathcal{G}_{\cB_i}(v,\bj) &= \frac{1}{2}
  \int_{\theta_i-\alpha_i}^{\theta_i + \alpha_i} d \theta
  \int_{L-d(\theta)}^L \frac{dr}{r} \left\{ \int_r^L \frac{ds}{s}
\left[ \cos (k \theta) \partial_\theta^2 w_k(s,\theta) - \right. \right. \nonumber \\
 & \qquad \qquad \left. \left. 
2 k \sin(k \theta) \partial_\theta w_k(s,\theta) + \cU_i \partial_\theta^2 w(s,\theta)
\right]\right\}^2\,  \nonumber \\
&\le \frac{3}{2} \left[ \mathcal{S}_{i,1} + \mathcal{S}_{i,2} + 
\mathcal{S}_{i,3}\right]\, .
\end{align}
Here we let 
\begin{equation}
w_k(s,\theta) = \frac{(s/L)^k - [1-d(\theta)/L]^{2k} (L/s)^k}{1 
- [1-d(\theta)/L]^{2k}}\, ,\quad
w(s,\theta) = \frac{\ln(s/L)}{\ln[1-d(\theta)/L]}\, , \label{eq:w}
\end{equation}
used the inequality 
\[
(a+b+c)^2 \le 3(a^2 + b^2 + c^2)\, , \quad \forall a,b,c \in
\mathbb{R}\, ,
\]
and introduced the integrals
\begin{align}
\mathcal{S}_{i,1} &= \int_{\theta_i-\alpha_i}^{\theta_i + \alpha_i}
  d \theta \int_{L-d(\theta)}^L \frac{dr}{r} \left[ 2 k \sin(k
    \theta)\int_r^L \frac{ds}{s}\, 
    \partial_\theta w_k(s,\theta)\right]^2, \label{eq:S1}\\
  \mathcal{S}_{i,2} &= \int_{\theta_i-\alpha_i}^{\theta_i + \alpha_i}
  d \theta \int_{L-d(\theta)}^L \frac{dr}{r} \left[ \cos (k \theta)
    \int_r^L \frac{ds}{s}\,
    \partial_\theta^2 w_k(s,\theta)\right]^2, \label{eq:S2}\\
\mathcal{S}_{i,3} &= 
  \int_{\theta_i-\alpha_i}^{\theta_i + \alpha_i} d \theta
  \int_{L-d(\theta)}^L \frac{dr}{r} \left[\cU_i \int_r^L \frac{ds}{s}
\,    \partial_\theta^2 w(s,\theta)\right]^2. \label{eq:S3}\\
\end{align} 

\subsection{Estimate of (\ref{eq:S1})}
We obtain from definition (\ref{eq:w}) that
\begin{equation}
  \partial_\theta w_k(s,\theta) = -\frac{2 k d'(\theta) 
p^2(\theta)}{(L-d(\theta)) (1-p^2(\theta))^2} \left[\left(\frac{s}{L}
\right)^k-\left(\frac{L}{s}\right)^k\right]\, ,
\label{eq:dwk}
\end{equation}
with 
\begin{equation}
\label{eq:defp}
p(\theta) = [1-d(\theta)/L]^k\, ,
\end{equation}
so we can bound $\mathcal{S}_{i,1}$ as 
\begin{align*}
  \mathcal{S}_{i,1} &\le 16 \int_{\theta_i-\alpha_i}^{\theta_i +
    \alpha_i} \hspace{-0.05in}d \theta \, \left[ \frac{k
      d'(\theta)p^2(\theta)}{(L-d(\theta))(1-p^2(\theta))^2} \right]^2
  \hspace{-0.05in}\int_{L-d(\theta)}^L\hspace{-0.05in} \frac{dr}{r}
  \left[ \int_{r}^L \hspace{-0.05in}ds \, \left( \frac{k s^{k-1}}{L^k}
      - \frac{k L^k}{s^{k+1}}\right) \right]^2.
\end{align*}

The integral in $r$ is estimated by 
\begin{align}
  \int_{L-d(\theta)}^L \frac{dr}{r} \left[ \int_{r}^L
    \hspace{-0.05in}ds \, \left( \frac{k s^{k-1}}{L^k} - \frac{k
        L^k}{s^{k+1}}\right) \right]^2 = \int_{L-d(\theta)}^L
  \frac{dr}{r} \left[2 - \frac{r^k}{L^k} -
    \frac{L^k}{r^k}\right]^2 \nonumber \\
  \le - \frac{(1-p(\theta))^4}{p^2(\theta)} \ln \left[
      1-\frac{d(\theta)}{L} \right] , \label{eq:intR}
\end{align}
where  the monotonicity in $r$ of the function in parenthesis implies
\begin{align}
  \left(2 - \frac{r^k}{L^2} - \frac{L^k}{r^k}\right)^2 &\le \left[2 -
    \frac{(L-d(\theta))^k}{L^2} - \frac{L^k}{(L-d(\theta))^k}\right)^2
  = \frac{(1-p(\theta))^4}{p^2(\theta)}\,,
\label{eq:ineqQ}
\end{align}
for all $r \in [L-d(\theta),L]$.   Moreover, since
\begin{equation}
\label{eq:ineqQ1}
\frac{1}{1-p^2(\theta)} \le \frac{1}{1-p(\theta)}\, ,
\end{equation}
we obtain the bound
\begin{align}
  \mathcal{S}_{i,1} \le 16 \int_{\theta_i-\alpha_i}^{\theta_i +
    \alpha_i} d \theta \,
  \left[\frac{d'(\theta)}{L-d(\theta)}\right]^2 \left\{-\left[k
    p(\theta)\right]^2\ln\left[ 1-\frac{d(\theta)}{L} \right] \right\}\, .
  \label{eq:Si1}
\end{align}
Function $k p(\theta)$ attains its maximum at $k = -1/\ln\left[
  1-{d(\theta)}/{L} \right]$
\begin{equation}
\label{eq:ineqQ2}
k p(\theta) = k \left[1-\frac{d(\theta)}{L}\right]^k \le \frac{e^{-1}}{-\ln\left[
  1-\frac{d(\theta)}{L} \right]}\,, 
\end{equation}
and after expanding the logarithm we get
\begin{equation}
\label{eq:EstS1_1}
\mathcal{S}_{i,1} \le C \int_{\theta_i-\alpha_i}^{\theta_i +
    \alpha_i} d \theta \,
  \frac{[d'(\theta)]^2}{L d(\theta)} \, ,
\end{equation}
with  positive constant $C$ of order one.

To estimate (\ref{eq:EstS1_1}) we obtain from definition 
(\ref{eq:dTh}) that 
\begin{align}
d(\theta) &= L -\rho_i \cos(\theta-\theta_i) - \sqrt{R^2 -
\rho_i^2 \sin^2(\theta-\theta_i)} \nonumber \\
&= \de_i + \rho_i[1 - \cos(\theta-\theta_i)] + 
R - \sqrt{R^2 -
\rho_i^2 \sin^2(\theta-\theta_i)} \nonumber \\
&\ge R - \sqrt{R^2 -
\rho_i^2 \sin^2(\theta-\theta_i)}\, ,
\label{eq:bound1}
\end{align}
and note that its derivative satisfies 
\begin{equation}
  |d'(\theta)| = \frac{\rho_i \sin |\theta-\theta_i| }{ \sqrt{R^2 -
      \rho_i^2 \sin^2(\theta-\theta_i)}} \left[ L - d(\theta)\right] 
  \le \frac{2L}{R} \rho_i \left|\sin
    (\theta-\theta_i)\right|\, .
\label{eq:bound3}
\end{equation}
Here we used (\ref{eq:alpha}) to write
\[
\frac{1}{\sqrt{R^2-\rho_i^2 \sin^2(\theta - \theta_i)}} \le
\frac{2}{R}, \quad \forall \, \theta \in
(\theta_i-\alpha_i,\theta_i+\alpha_i)\,.
\]
The second derivative of $d(\theta)$, needed in the next
section, is bounded similarly
\begin{equation}
|d''(\theta)| \le \frac{8 L^2}{R} \, .
\label{eq:d2b}
\end{equation}

Inequalities (\ref{eq:bound1})-(\ref{eq:bound3}) give
\begin{align}
  \frac{\left[d'(\theta)\right]^2}{L d(\theta)} &\le \frac{4 L}{R^2}
  \frac{\rho_i^2 \sin^2(\theta-\theta_1)}{R - \sqrt{R^2 - \rho_i^2
      \sin^2(\theta-\theta_i)}} \nonumber \\
&= \frac{4 L}{R^2} \left[R + \sqrt{R^2 - \rho_i^2
      \sin^2(\theta-\theta_i)}\right] \nonumber \\
&=\left(\frac{8 L}{R}\right)\, ,
\label{eq:bound4}
\end{align}
and the estimate 
\begin{equation}
\mathcal{S}_{i,1} \le O(1) 
\label{eq:EstS1}
\end{equation}
follows from (\ref{eq:EstS1_1}) and $\alpha_i = O(R/L)$. 

\subsection{Estimate of (\ref{eq:S2})}
We obtain from (\ref{eq:dwk}) that 
\begin{align*}
  \partial_\theta^2 w_k(s,\theta) = \left\{ \frac{2[d'(\theta)
      p(\theta)]^2\left[ (2k+1)p^2(\theta)
        +2k-1\right]}{[1-p^2(\theta)]^3 [L-d(\theta)]^2} - \frac{2
      d''(\theta)p^2(\theta)}{ [1-p^2(\theta)]^2
      [L-d(\theta)]}\right\} \times \\
  \left[ \frac{k s^k}{L^k} - \frac{kL^k}{s^k}\right]\, ,
\end{align*}
and using the estimate (\ref{eq:intR}) of the integral in $r$ we get
\begin{align*}
  \mathcal{S}_{i,2} \le C \int_{\theta_i-\alpha_i}^{\theta_i +
    \alpha_i} d \theta \, \left\{- \frac{(1-p(\theta))^4}{p^2(\theta)}
    \ln \left[ 1-\frac{d(\theta)}{L} \right]\right\} \left\{ \left[
      \frac{d''(\theta)p^2(\theta)}{ [1-p^2(\theta)]^2
        [L-d(\theta)]}\right]^2 + \right. \nonumber \\
  \left. \left[\frac{[d'(\theta) p(\theta)]^2\left[
          (2k+1)p^2(\theta) +2k-1\right]}{[1-p^2(\theta)]^3
        [L-d(\theta)]^2}\right]^2 \right\}\, ,
\end{align*}
with positive constant $C$ of order one.  Now use inequality
(\ref{eq:ineqQ1}) and expand the logarithm and the terms
$L-d(\theta)$ in the denominator to simplify the bound
\begin{align*}
  \mathcal{S}_{i,2} \lesssim C \int_{\theta_i-\alpha_i}^{\theta_i +
    \alpha_i} d \theta \, \frac{d(\theta)}{L} \left\{ \left[
      \frac{d''(\theta)}{L}\right]^2 p^2(\theta) + 16 \left[
      \frac{\left[ d'(\theta) \right]^2}{L} \right]^2 \left[ \frac{k
        p(\theta)}{L[1-p(\theta)]}\right]^2 \right\}\, .
\end{align*}
The derivatives of $d(\theta)$ are estimated in
(\ref{eq:d2b}) and (\ref{eq:bound4}), $p(\theta)\le 1$,  and
\begin{equation}
\label{eq:monot}
\frac{k p(\theta)}{1-p(\theta)} = \frac{k \left[1-d(\theta)/L \right]^k}{
  1-\left[1-d(\theta)/L\right]^k} \le \frac{L}{d(\theta)}
\left[ 1-d(\theta)/L \right] \, .
\end{equation}
This is because the function $k a^{k}/(1-a^{k})$ is monotonically
decreasing in $k$ for any $a \in (0,1)$ and $k \ge 1$. In particular,
for $a = 1-d/L$ we have (\ref{eq:monot}).  Gathering all the results
and using that $\alpha_i = O(R/L)$ we get
\begin{equation}
  \mathcal{S}_{i,2} \le C_1 \int_{\theta_i-\alpha_i}^{\theta_i +
    \alpha_i} d \theta \, \frac{d(\theta)L}{R^2} \le O(1)\, .
\label{eq:EstS2}
\end{equation}
\subsection{Estimate of (\ref{eq:S3})}
We recall from (\ref{eq:MPR}) that $\cU_i$ is at most of order one, 
and obtain from (\ref{eq:w}) that 
\begin{align}
  \partial_\theta^2 w(s,\theta) &= \left\{ \frac{2
      [d'(\theta)]^2/\ln[1-d(\theta)/L] + [d'(\theta)]^2 +
      [L-d(\theta)] d''(\theta)}{[L-d(\theta)]^2 \left[
        \ln[1-d(\theta)/L] \right]^2} \right\} \ln
  \frac{s}{L}\,\nonumber \\
  &\approx \left\{ \frac{d''(\theta)}{L} - 2
    \frac{[d'(\theta)]}{d(\theta)L} \right\} \frac{\ln (s/L)}{ \left[
      \ln[1-d(\theta)/L] \right]^2}\, .
\end{align}
The integrals in $s$ and $r$ give  
\[
  \int_{L-d(\theta)}^L \frac{dr}{r} \left[ \int_r^L \frac{ds}{s}
    \,    \ln \left(\frac{s}{L}\right)\right]^2 = 
  \frac{1}{4} \int_{L-d(\theta)}^L \frac{dr}{r} \left[ \ln\frac{r}{L}\right]^4=
  -\frac{1}{20} 
\left\{ \ln\left[1 -\frac{d(\theta)}{L}\right] \right\}^5, 
\]
and with the bounds (\ref{eq:d2b}) and (\ref{eq:bound4}) of the
derivatives of $d(\theta)$, and the expansion of the logarithm, we
obtain the estimate
\begin{equation}
  \mathcal{S}_{i,3}\le \frac{9}{20} \int_{\theta_i-\alpha_i}^{
\theta_i+\alpha_i} d\theta \, \frac{L \, d(\theta) }{R^2} \le 
  O(1)\, .
\end{equation}
\section{Energy in the boundary layer}
\label{ap:BL_EN}
We use the test potential (\ref{eq:testV}) to calculate the upper 
bound of the energy in the boundary layer. Given the decomposition
of the layer in the sets $\cB_i$ and $\cB_{i+}$, we write the bound as in 
(\ref{eq:UB1}), and estimate the two terms in sections \ref{sect:Bi}
and \ref{sect:Bi+}.

\subsection{Energy in the sets $\cB_{i+}$}
\label{sect:Bi+}
Let us introduce the simplifying notation 
\begin{equation}
  \overline{\cU}_i = \frac{\cU_i+\cU_{i+1}}{2}, \qquad 
  \widetilde{\cU}_i = \cU_{i+1}-\cU_i\, , 
\label{eq:EB1}
\end{equation}
and 
\begin{equation}
  \overline{\theta}_i = \frac{\theta_i + \theta_{i+1}}{2} - 
  \frac{\alpha_{i+1} - \alpha_i}{2}\, , \qquad 
  \widetilde \theta_i = 
  (\theta_{i+1}-\alpha_{i+1}) - (\theta_i +\alpha_i)\,, 
\end{equation}
for the average and difference potentials and angles, so that 
\begin{equation}
\mathfrak{L}(\theta,\bU^\Gamma) = \overline{\cU}_i + \widetilde{\cU}_i 
\left(\frac{\theta-\overline{\theta}_i}{\widetilde \theta_i} \right)
\end{equation}
in $\cB_{i+}$. We obtain after straightforward calculation that
\begin{equation}
  \overline{E}_{\cB_{i+}}(\bU^\Gamma,\psi;v) = \frac{k \widetilde
    \theta_i}{4} + \mathcal{P}_{\cB_{i+}}\, ,
\label{eq:EB2}
\end{equation}
with perturbation term 
\begin{align}
\mathcal{P}_{\cB_{i+}} = &
\frac{k p^2 \widetilde \theta_i }{4(1-p^2)} -
  \frac{k^2 p^2 \ln\left(1-\frac{R}{2L}\right)} {2(1-p^2)^2}
  \int_{-\widetilde \theta_i/2}^{\widetilde \theta_i/2} d \theta
  \cos\left[2 k (\overline \theta_i + \theta)\right]+ \nonumber \\
  & \frac{1}{2\ln\left(1-\frac{R}{2L}\right)} \int_{-\widetilde
    \theta_i/2}^{\widetilde \theta_i/2} d \theta
  \left(\overline{\cU}_i + \widetilde{\cU}_i \frac{\theta}{\widetilde
      \theta_i} \right) \cos\left[k(\overline \theta_i +
    \theta)\right] - \nonumber \\
  & \frac{1}{4\ln\left(1-{R}/{2L}\right)} \int_{-\widetilde
    \theta_i/2}^{\widetilde \theta_i/2} d \theta
  \left(\overline{\cU}_i + \widetilde{\cU}_i \frac{\theta}{\widetilde
      \theta_i} \right)^2 -\frac{\widetilde \cU_i^2
    \ln\left(1-\frac{R}{2L}\right)}{6
    \widetilde \theta_i} + \nonumber \\
  &\frac{\widetilde \cU_i \left[(1-p^2)+ 2 p
      \ln p\right]}{2(1-p^2)
    \ln p}\, \frac{1}{\widetilde  \theta_i}
  \int_{-\widetilde \theta_i/2}^{\widetilde \theta_i/2} d \theta
  \sin\left[k(\overline \theta_i + \theta)\right]\, ,
\label{eq:EB3}
\end{align}
where
\[
p = \left[1-\frac{d(\theta)}{L}\right]^k =
\left(1-\frac{R}{2L}\right)^k\, .
\]

Now let us show that $\mathcal{P}_{\cB_{i+}} = O(1)$. The first term
in (\ref{eq:EB3}) is estimated as
\[
\frac{k p^2 \widetilde \theta_i }{(1-p^2)} = \frac{k
  \left(1-\frac{R}{2L}\right)^{2k} \widetilde \theta_i }{
  \left[1-\left(1-\frac{R}{2L}\right)^{2k}\right]} \le
\frac{\left(1-\frac{R}{2L}\right)^{2} \widetilde \theta_i }{
  \left[1-\left(1-\frac{R}{2L}\right)^{2}\right]} = O(1)\, ,
\]
because the function is monotonically decreasing in $k$ and
$\widetilde \theta_i = O(R/L)$.  The second term in (\ref{eq:EB3})
satisfies
\begin{align*}
  \left| \frac{k^2 p^2 \ln\left(1-\frac{R}{2L}\right)} {(1-p^2)^2}
    \int_{-\widetilde \theta_i/2}^{\widetilde \theta_i/2} d \theta
    \cos\left[2 k (\overline \theta_i + \theta)\right]\right| = \left(
    \frac{kp}{1-p^2}\right)^2 \left|\cos(k \overline \theta_i)
    \mbox{sinc}\left(\frac{k \widetilde \theta_i}{2}\right)\right|
  \times \\ \widetilde \theta_i \left|
    \ln\left(1-\frac{R}{2L}\right)\right| \le \left(
    \frac{kp}{1-p}\right)^2 \widetilde \theta_i \left|
    \ln\left(1-\frac{R}{2L}\right)\right| \le O(1)\, ,
\end{align*}
where we used the inequality (\ref{eq:monot}) and expanded the
logarithm. The third, fourth and fifth terms in (\ref{eq:EB3}) are
also order one, because the integrands are order one and
\[
\widetilde \theta_i \sim -\ln\left(1-\frac{R}{2L}\right) = O
\left(\frac{R}{L}\right)\, .
\]
The last term in (\ref{eq:EB3}) satisfies
\[
  \left|\frac{\widetilde \cU_i \left[(1-p^2)+ 2 p
        \ln p\right]}{(1-p^2)
      \ln p}\right|\, \frac{1}{\widetilde  \theta_i}
  \int_{-\widetilde \theta_i/2}^{\widetilde \theta_i/2} d \theta
  \sin\left[k(\overline \theta_i + \theta)\right] \le 
  \left|\frac{ \left[(1-p^2)+ 2 p
        \ln p\right]}{(1-p^2)
      \ln p}\right| \le 1\, ,
\]
because the potentials satisfy the maximum principle (\ref{eq:MPR}).
The last inequality is easy to see, for example by plotting the
function for $p \in (0,1)$. 
\subsection{Energy in the sets $\cB_{i}$}
\label{sect:Bi}
The test potential in this set is of the form 
\begin{equation}
v(r,\theta) = w_k(r,\theta) \cos (k \theta) + \cU_i w(r,\theta)\, , 
\end{equation}
with functions $w_k(r,\theta)$ and $w(r,\theta)$ defined in
(\ref{eq:w}). We write the contribution of $\cB_i$ to the energy bound
as a quadratic polynomial in the potentials
\begin{equation}
\overline{E}_{\cB_{i}}(\bU^\Gamma,\psi;v) = 
\frac{1}{2} \int_{\cB_i} d\bx \, \left| \nabla v(\bx)\right|^2 = 
a_i \cU_i^2 + 2 b_i \cU_i + c_i\, . 
\label{eq:E1}
\end{equation}
The leading coefficients are independent of $k$
\begin{equation}
  a_i = \frac{1}{2}  \int_{\cB_i} d\bx \,\left\{ 
    \left[\partial_r w(r,\theta)\right]^2 + 
    \left[\frac{1}{r} \partial_\theta w(r,\theta)\right]^2 \right\}\, ,
  \label{eq:defa} 
\end{equation}
and are estimated in section \ref{sect:esta}.
The coefficients of the linear term are
\begin{align}
  b_i = \frac{1}{2} \int_{\cB_i} d\bx \,\left\{ \cos (k \theta)
    \partial_r w_k(r,\theta)\partial_r w(r,\theta) - \frac{1}{r}
    \partial_\theta w(r,\theta)\times \right. \nonumber \\
  \left. \left[ \frac{k\sin (k \theta)}{r}  w_k(r,\theta) -
      \frac{ \cos (k \theta)}{r} \partial_\theta w_k(r,\theta)
    \right]\right\}\, ,
\label{eq:defb} 
\end{align}
and are estimated in section \ref{sect:estb}. The coefficients 
\begin{align}
  c_i = \frac{1}{2} \int_{\cB_i} \hspace{-0.05in}d\bx \,\left\{
    \left[\partial_r w_k(r,\theta)\right]^2 \cos^2(k\theta)+
    \left[\frac{k \sin(k \theta)}{r} w_k(r,\theta) - \frac{\cos (k
        \theta)}{r}
      \partial_\theta w_k(r,\theta)\right]^2 \right\}\, 
\label{eq:defc} 
\end{align}
are estimated in section \ref{sect:estc}.
\subsubsection{Estimate of $a_i$}
\label{sect:esta}
We obtain from (\ref{eq:w}) and (\ref{eq:defa}) after integrating in
the radial direction that
\begin{equation}
  a_i = \frac{1}{2} \int_{\theta_i-\alpha_i}^{\theta_i+\alpha_i} 
  \frac{d \theta}{-\ln \left[1-d(\theta)/L\right]} +
  \mathcal{P}_{B_i,a_i} 
  \, ,
\label{eq:A1}
\end{equation}
with remainder 
\begin{equation}
  \mathcal{P}_{B_i,a_i} = -\frac{1}{2} 
\int_{\theta_i-\alpha_i}^{\theta_i+\alpha_i} 
    d \theta \, \frac{[d'(\theta)]^2}{3 [L-d(\theta)]^2 \ln\left[
1-d(\theta)/L\right]}\, .
\end{equation}
We can bound it as
\begin{equation}
\left|\mathcal{P}_{B_i,a_i}\right| \le O(1)\, ,
\end{equation}
using the estimate (\ref{eq:bound4}) of $d'(\theta)$,  expanding
the logarithm and recalling that the angle $\alpha_i = O(R/L)$. 

To calculate the first term in (\ref{eq:A1}), we expand the logarithm
\begin{align}
\label{eq:E4}
  \int_{\theta_i-\alpha_i}^{\theta_i+\alpha_i} \hspace{-0.05in}
  \frac{d \theta}{-\ln \left[1-d(\theta)/L\right]} &=L \left[ 1+
    o(1)\right] \int_{-\alpha_i}^{\alpha_i} \frac{d \theta
  }{d(\theta_i + \theta)} \, ,
\end{align}
and obtain an integral that is basically the same as that in
(\ref{eq:OLDINT}). Recalling definition (\ref{eq:bound1}) of
$d(\theta)$ and using that $\alpha_i = O(R/L)$, we have the
approximation
\begin{equation}
  \frac{L}{d(\theta_i + \theta)} = \frac{L}{\de_i + 
\frac{\rho_i L}{2 R} \theta^2} + 
  O\left(\frac{L}{R}\right)\, ,
\label{eq:approxd}
\end{equation}
and the coefficient becomes
\begin{align}
  a_i &= \frac{1}{2} \int_{-\alpha_i}^{\alpha_i} d \theta \, 
\frac{L}{\de_i + \frac{\rho_i L}{2 R} \theta^2} + O(1) \nonumber \\
  &= \frac{1}{2} \sqrt{\frac{2 L R}{\rho_i \de_i}} \int_{-\alpha_i
    \sqrt{\frac{\rho_i L}{2 R \de_i}}}^{\alpha_i \sqrt{\frac{\rho_i
        L}{2 R \de_i}}} \frac{dt}{1+t^2} + O(1) \, ,
  \nonumber \\
  &= \frac{\pi}{2} \sqrt{\frac{\rho_i L}{2 R \de_i}} + O(1)\, .
\label{eq:a1}
\end{align}

\subsubsection{Estimate of $b_i$}
\label{sect:estb}
We obtain from (\ref{eq:w}) and (\ref{eq:defb}) after integrating in 
the radius that 
\begin{align}
  b_i = \frac{1}{2}\int_{\theta_i - \alpha_i}^{\theta_i + \alpha_i}
  \hspace{-0.05in} d \theta \left\{ \frac{\cos (k \theta)}{\ln
      \left[1-d(\theta)/L\right]} + \frac{k d'(\theta) \sin(k
      \theta)}{L-d(\theta)]} \, \frac{[1-p^2(\theta) + 2 p(\theta)\ln
      p(\theta)]}{[1-p^2(\theta)][\ln p(\theta)]^2}+
  \right. \nonumber \\
  \left. - \frac{2 [d'(\theta)]^2 \cos (k \theta)}{[L-d(\theta)]^2
      \ln[ 1- d(\theta)/L]}\, \frac{p(\theta)\left[ 1 - p^2(\theta) +
        (1+p^2(\theta)) \ln p(\theta)\right]}{[1-p^2(\theta)]^2 \ln
      p(\theta)} \right\}\, .
\label{eq:bb1}
\end{align}
We show next that the first term may be large, but the last two 
are at most order one.

We estimate the last term in (\ref{eq:bb1}) using (\ref{eq:bound4})
and the bound
\begin{align*}
  \frac{p\left[ 1 - p^2 + (1+p^2) \ln p\right]}{(1-p^2)^2 \ln p} &\le
  \lim_{p\to 1} \frac{p\left[ 1 - p^2 + (1+p^2) \ln
      p\right]}{(1-p^2)^2 \ln p} = \frac{1}{6}\, ,
\end{align*}
which holds because the function on the left is monotonically
increasing in the interval $p \in (0,1)$.  We have 
\begin{align}
  \left|\int_{\theta_i - \alpha_i}^{\theta_i + \alpha_i}
    \hspace{-0.05in}d \theta \, \frac{2 [d'(\theta)]^2 \cos (k
      \theta)}{[L-d(\theta)]^2 \ln[ 1- d(\theta)/L]}\,
    \frac{p(\theta)\left[ 1 - p^2(\theta) + (1+p^2(\theta)) \ln
        p(\theta)\right]}{[1-p^2(\theta)]^2 \ln
      p(\theta)}  \right| \le \nonumber \\
  \frac{1}{6} \int_{\theta_i - \alpha_i}^{\theta_i + \alpha_i} d
  \theta \, \frac{[d'(\theta)]^2}{L d(\theta) [1-d(\theta)/L]^2} \le
  O(1)\, ,
\end{align} 
where we expanded the logarithm, and used that $\alpha_i = O(R/L)$.

The second term in (\ref{eq:bb1}) is estimated using integration by
parts
\begin{align*}
  &\int_{\theta_i - \alpha_i}^{\theta_i + \alpha_i} d \theta \,
  \frac{k d'(\theta) \sin(k \theta)}{L-d(\theta)]} \,
  \frac{[1-p^2(\theta) + 2 p(\theta)\ln
    p(\theta)]}{[1-p^2(\theta)][\ln p(\theta)]^2} =   \\
  & \qquad \left. -\frac{d'(\theta) \cos(k \theta)}{L-d(\theta)]} \,
    \frac{[1-p^2(\theta) + 2 p(\theta)\ln p(\theta)]}{[1-p^2(\theta)]
      [\ln p(\theta)]^2}\right|_{\theta_{i}-\alpha_i}^{\theta_{i}+
    \alpha_i} - \\
  & \qquad \int_{\theta_i - \alpha_i}^{\theta_i + \alpha_i} d \theta
  \, \cos(k \theta) \frac{\left[[L-d(\theta)]d''(\theta) +
      [d'(\theta)]^2\right]}{[L-d(\theta)]^2} \, \frac{[1-p^2(\theta)
    + 2
    p(\theta)\ln p(\theta)]}{[1-p^2(\theta)][\ln p(\theta)]^2} - \\
  & \qquad \int_{\theta_i - \alpha_i}^{\theta_i + \alpha_i} d \theta
  \, \cos(k \theta) \frac{d'(\theta)}{[L-d(\theta)]} \,
  \frac{d}{d\theta}\left\{\frac{[1-p^2(\theta) + 2 p(\theta)\ln
      p(\theta)]}{[1-p^2(\theta)][\ln p(\theta)]^2}\right\}\, .
\end{align*}
We have that 
\begin{align*}
  \frac{\left[ 1 - p^2 + 2 \ln p\right]}{(1-p^2)^2 (\ln p)^2} &\le
  \lim_{p\to 1} \frac{\left[ 1 - p^2 + 2 \ln
      p\right]}{(1-p^2)^2 (\ln p)^2} = \frac{1}{6}\, ,
\end{align*}
because the function is monotonically increasing in the interval
$p \in (0,1)$. Moreover, 
\begin{align*}
  &\frac{d}{d\theta}\left\{\frac{[1-p^2(\theta) + 2 p(\theta)\ln
      p(\theta)]}{[1-p^2(\theta)][\ln p(\theta)]^2} \right\}= \frac{2
    d'(\theta)}{[L-d(\theta)] \ln[1-d(\theta)/L]} \times \\
  & \qquad \frac{[1-p^2(\theta)]^2 + p(\theta)[1-p^2(\theta)] 
\ln p(\theta)
    -p(\theta) [1+p^2(\theta)][\ln p(\theta)]^2}{[1-p^2(\theta)]^2
    [\ln p(\theta)]^2}\, 
\end{align*}
with the last factor bounded in the interval $p \in (0,1)$, as can be
seen easily by plotting. Thus, gathering all the results and using 
the estimates (\ref{eq:d2b}) and (\ref{eq:bound4}) for the derivatives
of $d(\theta)$, we get 
\begin{align}
  \left|\int_{\theta_i - \alpha_i}^{\theta_i + \alpha_i} d \theta \,
  \frac{k d'(\theta) \sin(k \theta)}{L-d(\theta)]} \,
  \frac{[1-p^2(\theta) + 2 p(\theta)\ln
    p(\theta)]}{[1-p^2(\theta)][\ln p(\theta)]^2}\right| \le O(1)\, .
\end{align}

We have obtained that 
\begin{align}
  b_i &= \frac{1}{2}\int_{\theta_i- \alpha_i}^{\theta_i + \alpha_i} d
  \theta \,
  \frac{\cos (k \theta)}{\ln[1- d(\theta)/L]} + O(1) \nonumber \\
  &= - \frac{L}{2}\cos (k \theta_i) \int_{- \alpha_i}^{\alpha_i} d
  \theta \, \frac{\cos (k \theta)}{ d(\theta_i + \theta)} + O(1)
  \nonumber \\
  &= -\frac{L}{2 \de_i}\cos(k \theta_i) \int_{- \alpha_i}^{\alpha_i} d
  \theta \, \frac{\cos (k \theta)}{1 + \frac{\rho_i L}{2 R \de_i}
    \theta^2} + O(1)\, ,
\label{eq:bb2}
\end{align}
where we expanded the logarithm and the cosine, discarded the sin
term which is odd, and used approximation (\ref{eq:approxd}).  The
integral in (\ref{eq:bb2}) is similar to that in (\ref{eq:E4}) for
small $k$, but for large $k$ the result is smaller due to the
oscillatory cosine.  Explicitly, we have
\begin{align}
  b_i &= -\frac{\cos(k \theta_i)}{2}\sqrt{\frac{2 RL}{\rho_i \de_i}}
  \int_{-\alpha_i \sqrt{\frac{\rho_i L}{2 R \de_i}}}^ {\alpha_i
    \sqrt{\frac{\rho_i L}{2 R \de_i}}} dt \, \frac{\cos\left(k
      \sqrt{\frac{2 R
          \de_i}{\rho_i L}}\right)}{1 + t^2} + O(1) \nonumber \\
  &= -\frac{\pi\cos(k \theta_i)}{2}\sqrt{\frac{2 RL}{\rho_i \de_i}}
  e^{-k \sqrt{\frac{2 R \de_i}{\rho_i L}}} + O(1)\, .
\end{align}

\subsubsection{Estimate of $c_i$}
\label{sect:estc} 
We obtain from (\ref{eq:w}) and (\ref{eq:defc}) after integrating in
the radius that
\begin{align}
  c_i = &\frac{k \alpha_i}{2} +
  \frac{1}{2}\int_{\theta_i-\alpha_i}^{\theta_i+\alpha_i}
  \hspace{-0.05in} d \theta \, \left\{
    \frac{kp^2(\theta)}{[1-p^2(\theta)]} - \frac{2 k^2 p^2(\theta)
      \ln[1-d(\theta)/L] \cos(2 k
      \theta)}{[1-p^2(\theta)]^2} \right. + \nonumber\\
  & \frac{2 k d'(\theta) p(\theta)\sin (2 k \theta)}{[L-d(\theta)]} \,
  \frac{p(\theta)[1-p^2(\theta) + (1+p^2(\theta))\ln
    p(\theta)]}{[1-p^2(\theta)]^3} +  \nonumber \\
  &\left. \frac{2 k^2 [d'(\theta)]^2 p^2(\theta) cos^2(k \theta)
      \ln[1-d(\theta)/L]}{L - d(\theta)]^2 [1-p^2(\theta)]^2}\,
    \frac{[1-p^4(\theta) + 4 p^2(\theta) \ln
      p(\theta)}{[1-p^2(\theta)]^2 \ln p(\theta)}\right\}\, ,
\label{eq:cii}
\end{align}
and proceeding as in the previous two sections, we conclude that the last
two terms are at most order one. 

To calculate the first integral in (\ref{eq:cii}), let us introduce the 
notation 
\[
  p^2(\theta_i + \theta) = \left[1 - \frac{d(\theta_i
      +\theta)}{L}\right]^{2k} = e^{-2 k x}\, , \quad x = -\ln\left[1 -
    \frac{d(\theta_i +\theta)}{L}\right]\, , \quad \theta \in
  (-\alpha_i,\alpha_i)\, ,
\]
and use that $d \ll L$ to write 
\[
x = \frac{d}{L} + O\left(\frac{d^2}{L^2}\right) = 
\frac{\widetilde d}{L} + \frac{d - \widetilde{d}}{L} + 
O\left(\frac{d^2}{L^2}\right)\,,
\]
for 
\[
\widetilde d = \de_i + \frac{\rho_i L}{2 R} \theta^2\, ,
\]
the parabolic approximation of $d(\theta_i + \theta)$.
We have from (\ref{eq:approxd}) that
\[
|d - \widetilde d| \le O\left(\frac{d^2}{R}\right)\, ,
\]
and therefore 
\[
 x = \frac{\de_i}{L} \left[ 1 + \frac{\rho_i L}{2 R \de_i} \theta^2\right] 
+ O\left(\frac{d^2}{LR}\right)\, .
\]

Next, we let 
\[
\frac{k p^2}{1-p^2} = \frac{k e^{-2kx}}{1-e^{-2kx}} =
\frac{k}{e^{2kx}-1}=:f(k,x)\, ,
\]
and use the mean value theorem to write 
\[
f(k,x) = \frac{k}{e^{2 k \widetilde{d}/L} - 1} + \partial_x
f(k,x') \left[x-\frac{\widetilde d}{L}\right]\, , 
\]
for some $x' \sim \widetilde d/L$. Note that because $|\partial_x
f(k,x)|$ is monotonically decreasing in $k$, we have
\[
|\partial_k f(k,x)| = \frac{2 k^2 e^{2 k x}}{[e^{2kx}-1]^2}
\le |\partial_k f(x,k=1)| = \frac{2 e^{2 x}}{[e^{2x}-1]^2}
= O\left(\frac{1}{x^2}\right) = O\left(\frac{L^2}{d^2}\right)\,.
\]
Therefore, we can approximate 
\[
f(k,x) = \frac{k p^2}{1-p^2} = \frac{k}{e^{2 k \widetilde{d}/L} - 1} + 
O\left(\frac{L}{R}\right)\, ,
\]
and since $\alpha_i = O(R/L)$, we write the first integral in
(\ref{eq:cii}) as
\begin{align}
  \int_{-\alpha_i}^{\alpha_i} d \theta \, \frac{ k p^2(\theta_i +
    \theta)}{1-p^2(\theta_i + \theta)} &= \int_{-\alpha_i}^{\alpha_i}
  d \theta \, \frac{k}{e^{2 k\widetilde{d}/{L}} - 1} + O(1) \nonumber \\
  &= \frac{1}{2} \sqrt{\frac{2 R L}{\rho_i \de_i}} \int_{-Y_i}^{Y_i}
  dy \, \frac{\lambda}{e^{\lambda (1+y^2)} -1} + O(1)\, ,
  \label{eq:cii1}
\end{align}
with
\[
\lambda = \frac{2 k \de_i}{L}, \quad Y_i = \alpha_i \sqrt{\frac{ L
    \rho_i}{2 R \de_i}} \sim \sqrt{\frac{R}{\de_i}} \gg 1.
\]
Moreover,
\begin{align*}
  \int_{-Y_i}^{Y_i} dy \, \frac{\lambda}{e^{\lambda(1+y^2)} -1} &=
  \int_{-\infty}^{\infty} dy \, \frac{\lambda}{e^{\lambda (1+y^2)}
    -1} +
  O\left(\sqrt{\frac{\de}{R}} \right)\, \\
  & = \sqrt{\pi \lambda} \, \rm{Li}_{1/2}\left( e^{-\lambda} \right) +
  O\left(\sqrt{\frac{\de}{R}}\right)\, ,
\end{align*}
with remainder estimated as
\begin{align*}
  \int_{Y_i}^\infty dy \, \frac{\lambda}{e^{\lambda(1+y^2)} -1}
  \le \int_{Y_i}^\infty \frac{dy}{1+y^2} = \frac{\pi}{2} -
  \arctan(Y_i) = O\left(\frac{1}{Y_i}\right) = O\left(\sqrt{\frac{\de}{R}}\right)\,.
\end{align*}
Here we used that the integrand is monotonically decreasing in
$\lambda$ to write 
\[
\frac{\lambda}{e^{\lambda(1+y^2)} -1} \le \lim_{\lambda \to 0} 
\frac{\lambda}{e^{\lambda(1+y^2)} -1} = \frac{1}{1+y^2}\, .
\]

Gathering the results, we obtain that the first integral in
(\ref{eq:cii}) is given by
\begin{equation}
  \int_{\theta_i-\alpha_i}^{\theta_i+\alpha_i}
  \hspace{-0.05in} d \theta \,
  \frac{k p^2(\theta)}{1-p^2(\theta)} = \frac{\pi}{2} 
  \sqrt{\frac{2 R L}{\rho_i \de_i}}\, \sqrt{\frac{2 k \de_i}{L \pi}} 
  {\rm Li}_{1/2} \left( 
    e^{-2 k \de_i/L} \right) + O(1)\, .
\label{eq:int1}
\end{equation}
The second integral is obtained similarly, so we write directly its
expression
\begin{align}
  - \int_{\theta_i-\alpha_i}^{\theta_i+\alpha_i} \hspace{-0.05in} d
  \theta \, \frac{2 k^2 p^2(\theta) \ln[1-d(\theta)/L] \cos(2 k
    \theta)}{[1-p^2(\theta)]^2} = \pi \sqrt{\frac{2 R L}{\rho_i
      \de_i}} \cos^2(k \theta_i) e^{-2k \sqrt{\frac{2
        R \de_i}{\rho_i L}}} - \nonumber \\
  \frac{\pi}{2} \sqrt{\frac{2 R L}{\rho_i \de_i}} e^{-2 k
    \sqrt{\frac{2 \de_i R}{\rho_i L}}} + O(1)\, .
\label{eq:int2}
\end{align}
The result stated in (\ref{eq:UB2}) follows.

\newpage
\bibliographystyle{plain} \bibliography{BIBLIO}

\end{document}